\documentclass[10pt,journal]{IEEEtran}

\usepackage{xspace,amssymb,epsfig,subfigure,syntonly}
\usepackage{epsfig,amsmath,color}
\usepackage{xspace,syntonly}
\usepackage{url}
\usepackage{algorithm,algorithmic}

\newcommand{\utwi}[1]{\mbox{\boldmath $#1$}}

\newcommand{\trace}{{\textrm{Tr}}}
\newcommand{\rank}{{\textrm{rank}}}

\newcommand{\diag}{{\textrm{diag}}}

\newcommand{\cL}{{\cal{L}}}
\newcommand{\cN}{{\cal N}}

\newcommand{\cT}{{\cal T}}

\newcommand{\cE}{{\cal E}}

\newcommand{\cF}{{\cal F}}

\newcommand{\cU}{{\cal U}}

\newcommand{\cH}{{\cal H}}

\newcommand{\bc}{{\bf c}}

\newcommand{\bb}{{\bf b}}
\newcommand{\bd}{{\bf d}}
\newcommand{\be}{{\bf e}}

\newcommand{\bp}{{\bf p}}
\newcommand{\bq}{{\bf q}}

\newcommand{\bx}{{\bf x}}
\newcommand{\bu}{{\bf u}}
\newcommand{\bv}{{\bf v}}
\newcommand{\bi}{{\bf i}}

\newcommand{\bA}{{\bf A}}
\newcommand{\bB}{{\bf B}}
\newcommand{\bC}{{\bf C}}

\newcommand{\bL}{{\bf L}}
\newcommand{\bM}{{\bf M}}

\newcommand{\bS}{{\bf S}}
\newcommand{\bQ}{{\bf Q}}

\newcommand{\bI}{{\bf I}}

\newcommand{\bZ}{{\bf Z}}

\newcommand{\bY}{{\bf Y}}
\newcommand{\bV}{{\bf V}}

\newcommand{\bnu}{{\utwi{\nu}}}

\newcommand{\bPi}{{\utwi{\Pi}}}

\newcommand{\bmu}{{\utwi{\mu}}}

\usepackage{epstopdf}
\usepackage[noadjust]{cite}

\begin{document}

\newtheorem{definition}{Definition}
\newtheorem{remark}{Remark}
\newtheorem{proposition}{Proposition}
\newtheorem{lemma}{Lemma}

\title{Optimal Dispatch of Photovoltaic Inverters \\ in Residential Distribution Systems}
\author{Emiliano Dall'Anese, \emph{Member}, \emph{IEEE}, Sairaj V. Dhople, \emph{Member}, \emph{IEEE}, and Georgios B. Giannakis, \emph{Fellow}, \emph{IEEE}
\thanks{\protect\rule{0pt}{0.5cm}%
Submitted on July 14, 2013; revised on October 13, 2013 and November 21, 2013; accepted November 22, 2013.

This work was supported by the Institute of Renewable Energy and the Environment (IREE) grant no. RL-0010-13, University of Minnesota. 

The authors are with the Digital Technology Center and the Dept. of ECE, University of Minnesota, 200 Union Street SE, Minneapolis, MN 55455, USA. E-mails: {\tt \{emiliano, sdhople, georgios\}@umn.edu}
}
}

\markboth{IEEE TRANSACTIONS ON SUSTAINABLE ENERGY (TO APPEAR)}%
{Dall'Anese \MakeLowercase{\textit{et al.}}: }

\maketitle

\vspace{.5cm}

\begin{abstract}
Low-voltage distribution feeders were designed to sustain unidirectional power flows to residential neighborhoods. The increased penetration of roof-top photovoltaic (PV) systems has highlighted pressing needs to address power quality and reliability concerns, especially when PV generation exceeds the household demand. A systematic method for determining the active- and reactive-power set points for PV inverters in residential systems is proposed in this paper, with the objective of optimizing the operation of the distribution feeder and ensuring voltage regulation. Binary PV-inverter selection variables and nonlinear power-flow relations render the optimal inverter dispatch problem nonconvex and NP-hard. Nevertheless, sparsity-promoting regularization approaches and semidefinite relaxation techniques are leveraged to obtain a computationally feasible convex reformulation. The merits of the proposed approach are demonstrated using real-world PV-generation and load-profile data for an illustrative low-voltage residential distribution system.
\end{abstract}


\section{Introduction}
\label{sec:Introduction}

 \PARstart{T}{he} installed residential photovoltaic (PV) capacity increased by 61\% in 2012, driven in large by falling prices, increased consumer awareness, and governmental incentives~\cite{Sherwood12}. The proliferation of residential-scale roof-top PV systems presents a unique set of challenges related to power quality and reliability in low-voltage distribution systems. In particular, overvoltages experienced during periods when PV generation exceeds the household demand, and voltage sags during rapidly-varying irradiance conditions have become pressing concerns~\cite{Liu08,Tonkoski12}.

Efforts to ensure reliable operation of the existing distribution system with increased behind-the-meter PV generation are therefore focused on the possibility of PV inverters providing ancillary services~\cite{Turitsyn11}. This setup requires a departure from current interconnection standards~\cite{Standard1547}.
Organizations across Industry are endeavoring to address the issue and bring consistency to grid-interactive controls~\cite{Cpuc12,Nerc12}. For instance, by appropriately derating PV inverters, reactive power generation/consumption based on monitoring local electrical quantities has been recognized as a viable option to effect voltage regulation~\cite{Turitsyn11,Carvalho08,DeBrabandere04,Cagnano11,Farivar12,Aliprantis13}. However, such reactive power control (RPC) strategies typically yield low power factors (PFs) at the feeder input and high network currents, with the latter translating into power losses and possible conductor overheating~\cite{Tonkoski11}. To alleviate these concerns, an alternative is to curtail the active power produced by the PV inverters during peak irradiance hours~\cite{Tonkoski11,Tonkoski11Renewable}. For example, in~\cite{Tonkoski11}, the amount of active power injected by a PV inverter is lowered whenever its terminal voltage magnitude exceeds a predefined limit. The premise for these active power curtailment (APC) strategies is that the resistance-to-reactance ratios in low-voltage distribution networks renders the voltages very sensitive to variations in the active power injections~\cite{Tonkoski11,ADG-TPWRS12}. Of course, pertinent questions in this setup include what is the optimal amount of power to be curtailed, and by what PV systems in the network. 

A systematic and unified optimal inverter dispatch (OID) framework is proposed in this paper, with the goal of facilitating high PV penetration in existing distribution networks. The OID task involves solving an optimal power flow (OPF) problem to determine PV-inverter active- and reactive-power set points, so that the network operation is optimized according to well defined criteria (e.g., minimizing power losses), while ensuring voltage regulation and adhering to other electrical network constraints. The proposed OID framework provides increased flexibility over RPC or APC alone, by invoking a joint optimization of active and reactive powers. 
Indeed, while~\cite{DeBrabandere04,Carvalho08,Cagnano11,Turitsyn11,Aliprantis13,Tonkoski11} (and pertinent references therein) demonstrate the virtues of RPC and APC in effecting voltage regulation, these strategies are based on local information and, therefore, they may not offer system-level optimality guarantees. On the other hand, by promoting network-level optimization, the OPF-based formulation proposed in this work inherently ensures system-level benefits.

Unfortunately, due to binary PV-inverter selection variables and nonlinear power-balance constraints, the formulated problem is \emph{nonconvex} and \emph{NP-hard}. Nevertheless, a computationally-affordable \emph{convex} reformulation is derived here by leveraging sparsity-promoting regularization approaches~\cite{YuLi06,Wiesel11} and semidefinite relaxation techniques~\cite{Bai08, LavaeiLow,Tse12,Dallanese-TSG13}. Sparsity emerges because the proposed framework offers the flexibility of controlling only a (small) subset of PV inverters in the network. This allows one to discard binary optimization variables, and effect inverter selection by using group-Lasso-type regularization functions~\cite{YuLi06,Wiesel11}. 
To cope with nonconvexity due to bilinear power-balance constraints, the semidefinite relaxation approach~\cite{luospmag10} proposed in~\cite{Bai08, LavaeiLow} for the OPF problem in balanced transmission systems, and further extended to balanced and unbalanced distribution systems
in~\cite{Tse12, Lavaei_tree} and~\cite{Dallanese-TSG13}, respectively, is employed here too. This approach has the well documented merit of  finding the globally optimal solution of OPF problems in many practical setups~\cite{Tse12,Dallanese-TSG13}. The resultant relaxed OID is a semidefinite program (SDP), and it can be efficiently solved via primal-dual iterations or general purpose SDP solvers~\cite{Vandenberghe96}, with a significantly lower computational burden compared to off-the-shelf solvers for mixed-integer nonlinear programs (MINPs)~\cite{PaudyalyISGT}. Effectively, the OID task can be implemented in real
time, based on short-term forecasts of ambient conditions and loads. Subsequently, state-of-the-art advanced metering infrastructure can be used to relay set points to inverters~\cite{ADG-TPWRS12}. Towards this end, initiatives are currently underway to identify and standardize communication architectures for next-generation PV inverters~\cite{Key12}. 

To summarize, the main contributions of this work are as follows: i) A joint network-wide optimization of PV-inverter real and reactive power production improves upon conventional strategies that solely focus on either real~\cite{Tonkoski11} or reactive power~\cite{DeBrabandere04,Carvalho08,Cagnano11,Turitsyn11,Farivar12}. ii) Sparsity-promoting regularization approaches are proposed for the OPF problem to offer PV-inverter selection capabilities. While existing approaches either require controlling all the PV-inverters, or, assume that nodes providing ancillary services are preselected~\cite{DeBrabandere04,Carvalho08,Cagnano11,Turitsyn11,Tonkoski11,Farivar12,Dallanese-TSG13}, the formulated problem identifies optimal PV inverters for ancillary services provisioning. iii) A third contribution is represented by the RPC and APC strategies outlined in Section~\ref{sec:rpc}. In fact, the methods in~\cite{DeBrabandere04,Carvalho08,Cagnano11,Turitsyn11,Tonkoski11} are based on local information and assume that nodes providing ancillary services are preselected. In contrast, the RPC and APC tasks are cast here as instances of the OPF problem, and they offer inverter selection capabilities. OPF-type RPC strategies are proposed in~\cite{Farivar12,Bolognani13}; however, inverters providing reactive support are preselected. iv) The OID strategy is tested using real-world PV-generation and load-profile data for an illustrative low-voltage residential distribution system, and its performance is compared with RPC and APC.

The remainder of this manuscript is organized as follows. In Section~\ref{sec:Formulation}, the network and PV inverter models are described, and the RPC, APC, and OID strategies are introduced. Section~\ref{sec:Dispatch} outlines the proposed OID problem, while Section~\ref{sec:ConvexDispatch} describes its relaxed convex reformulation. Case studies and implementation details are presented in Section~\ref{sec:CaseStudies}, while Section~\ref{sec:Conclusions} concludes the paper and outlines future directions.

\section{Preliminaries and System Models}
\label{sec:Formulation}

\subsection{Notation}
Upper-case (lower-case) boldface
letters will be used for matrices (column vectors); $(\cdot)^\cT$ for transposition; $(\cdot)^*$ complex-conjugate; and, $(\cdot)^\cH$ complex-conjugate transposition;  
$\Re\{\cdot\}$ and $\Im\{\cdot\}$ denote the real and imaginary parts of a complex number, respectively; $j := \sqrt{-1}$ the imaginary unit. $\trace(\cdot)$ the matrix trace; $\rank(\cdot)$ the matrix rank; $|\cdot|$ denotes the magnitude of a number or the cardinality of a set; $\|\bv\|_2 := \sqrt{\bv^\cT \bv}$; $\|\bv\|_1 := \sum_i |[\bv]_i|$; and $\|\cdot\|_F$ stands for the Frobenius norm. Finally, $\bI_N$ denotes the $N \times N$ identity matrix; and, $\mathbf{0}_{M\times N}$, $\mathbf{1}_{M\times N}$ the $M \times N$ matrices with all zeroes and ones, respectively. 

\subsection{Distribution-Network Electrical Model}
\label{sec:DistributionNetworkModel}

\begin{figure}[t]
\begin{center}
\includegraphics[width=0.39\textwidth]{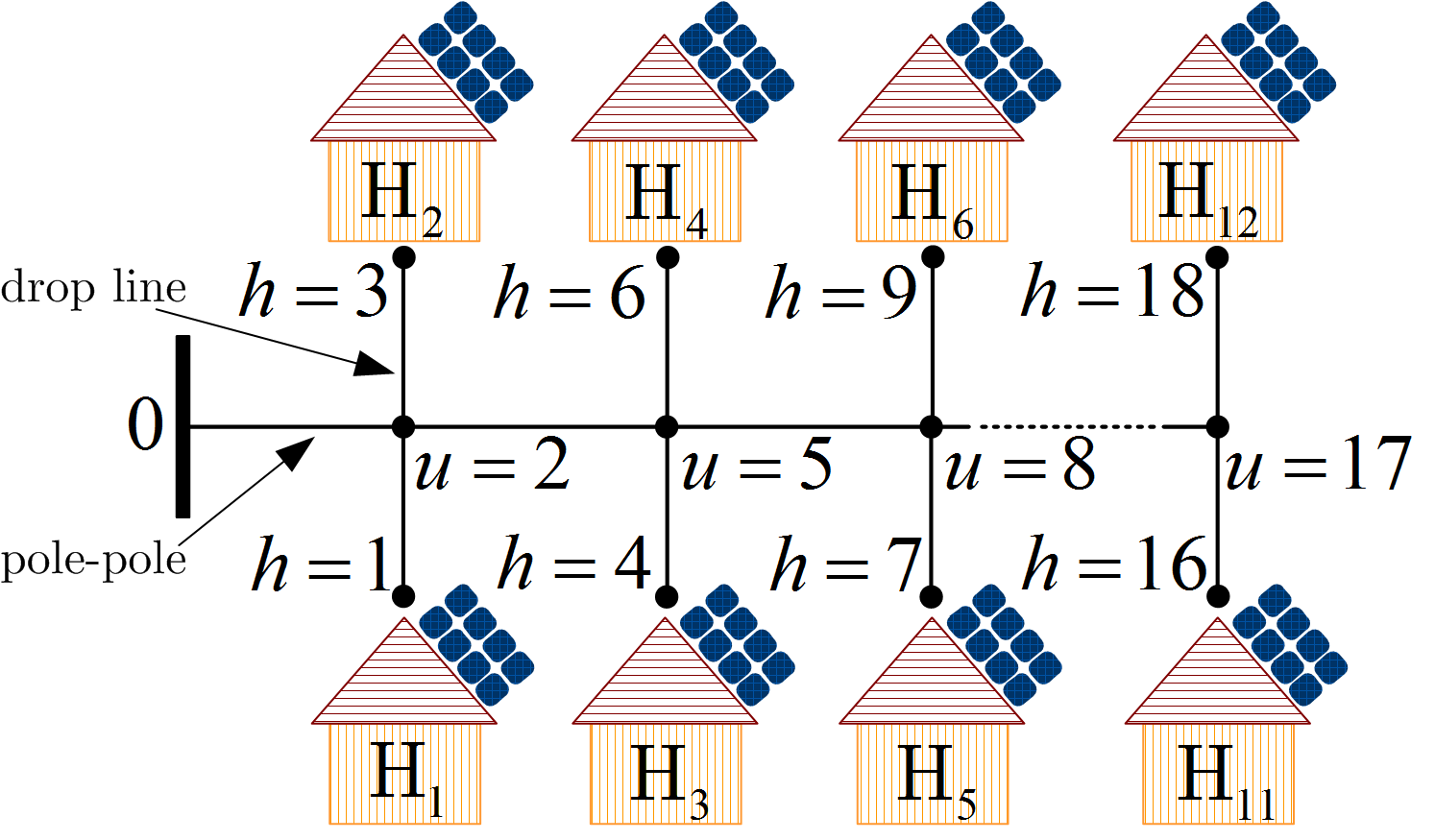}
\caption{Example of low-voltage residential network with high PV penetration adopted from~\cite{Tonkoski11, Tonkoski12}. Node 0 corresponds to the secondary of the step-down transformer, while set $\cU = \{2, 5, 8, 11, 14, 17\}$ collects nodes corresponding to distribution poles. Homes $\mathrm{H}_1,\dots,\mathrm{H}_{12}$ are connected to nodes in the set $\cH = \{1,3,4,6,7,9,10,12,13,15,16,18\}$.}
\label{F_LV_network}
\vspace{-.5cm}
\end{center}
\end{figure}

Consider a low-voltage radial feeder comprising $N+1$ nodes collected in the set
 $\cN := \{0,1,\ldots,N\}$, and overhead lines represented by the set of edges $\cE := \{(m,n)\} \subset \cN \times \cN$. Node $0$ is taken to be the secondary of the step-down transformer. Let $\cU, \cH \subset \cN$ collect nodes that correspond to utility poles, and houses with roof-top PV systems, respectively. For example, consider the residential network in Fig.~\ref{F_LV_network}, which is adopted from~\cite{Tonkoski11, Tonkoski12}, and is utilized in the case studies presented subsequently. In this network, one has that $\cU = \{2, 5, 8, 11, 14, 17\}$ and $\cH = \{1,3,4,6,7,9,10,12,13,15,16,18\}$ (the latter corresponding to houses $\mathrm{H}_1,\dots,\mathrm{H}_{12}$).

The distribution lines $(m,n) \in \cE$ are modeled as $\pi$-equivalent circuits. The series and shunt admittances of  line $(m,n)$ are given by $y_{mn} = (R_{mn} + j \omega L_{mn})^{-1}$ and $\bar{y}_{mn} = j \omega C_{mn}$, respectively, where $R_{mn}$, $L_{mn}$, $C_{mn}$, and $\omega$ denote the line resistance, inductance, shunt capacitance, and angular frequency~\cite{GSO08}. The distribution-network admittance matrix is denoted by $\bY \in \mathbb{C}^{N+1 \times N+1}$, and its entries are defined as 
\begin{equation}
[\bY]_{m,n} := \left \{\begin{array}{ll} 
\sum_{j \in \cN_m} \bar{y}_{mj} + y_{mj} , & \textrm{if } m=n\\ 
- y_{mn} , &\textrm{if }(m,n)\in \cE\\
0, & \textrm{otherwise} \end{array} \right.
\label{eq:Ymatrix}
\end{equation}
where $\cN_{m} := \{j \in \cN: (m,j) \in \cE\}$ denotes the set of nodes connected to the $m$-th one through a distribution line. Let $V_n, I_n \in \mathbb{C}$ denote the phasors for the line-to-ground voltage and the current injected at node $n \in \cN$, respectively. Collecting the currents injected in all nodes in the vector $\bi := [I_0, I_1, \ldots, I_N]^\cT \in \mathbb{C}^{N+1}$, and the node voltages in the vector $\bv := [V_0, V_1, \ldots, V_N]^\cT \in \mathbb{C}^{N+1}$, it follows that Ohm's law can be rewritten in matrix-vector form as 
\begin{equation}
\bi = \bY \bv .
\label{eq:OhmLaw} 
\end{equation}

\subsection{Residential PV Inverter and Load Models}
\label{sec:PVandLoadModel}
Given prevailing ambient conditions, let $\bar{P}_h$ denote the maximum active power that can be generated by the PV inverter(s) located at node $h \in \cH$. Henceforth, $\bar{P}_h$ is referred to as the \emph{available active power} at node $h$. The $|\cH| \times 1$ vector  collecting the available active powers $\{\bar{P}_h\}_{h \in \cH}$ is denoted by $\bar{\bp}_h$. Let $P_{s,h}$ and $Q_{s,h}$ denote the active power injected and the reactive power generated/consumed by the PV inverter at node $h$, respectively. The $|\cH| \times 1$ vectors collecting $\{P_{s,h}\}_{h \in \cH}$ and $\{Q_{s,h}\}_{h \in \cH}$ are denoted by $\bp_{s}$ and $\bq_s$, respectively. For conventional grid-tied residential-scale inverters that do not offer energy storage capabilities and operate at unity power factor, it follows that $P_{s,h} = \bar{P}_h$ and $Q_{s,h} = 0$. Nevertheless, since strategies where PV inverters are allowed to curtail their power output will be considered in the ensuing subsection, it is useful to denote the active power curtailed by inverter $h \in \cH$ as  
\begin{equation}
P_{c,h} := \bar{P}_h - P_{s,h} .
\label{eq:Power_curtailed}
\end{equation}
Furthermore, collect all the curtailed active powers in the $|\cH| \times 1$ non-negative vector $\bp_c$.  

A constant $PQ$ model is adopted for the residential loads, and $P_{\ell,h}$ and $Q_{\ell,h}$ denote the residential-load active and reactive powers at node $h \in \cH$. The corresponding $|\cH| \times 1$ vectors that collect the load active and reactive powers are denoted by $\bp_{\ell}$ and $\bq_{\ell}$, respectively. 

\subsection{PV-Inverter Control Strategies} 
\label{sec:Inverter}

\begin{figure}[t]
\centering
  \subfigure[]{\includegraphics[width=0.10\textwidth]{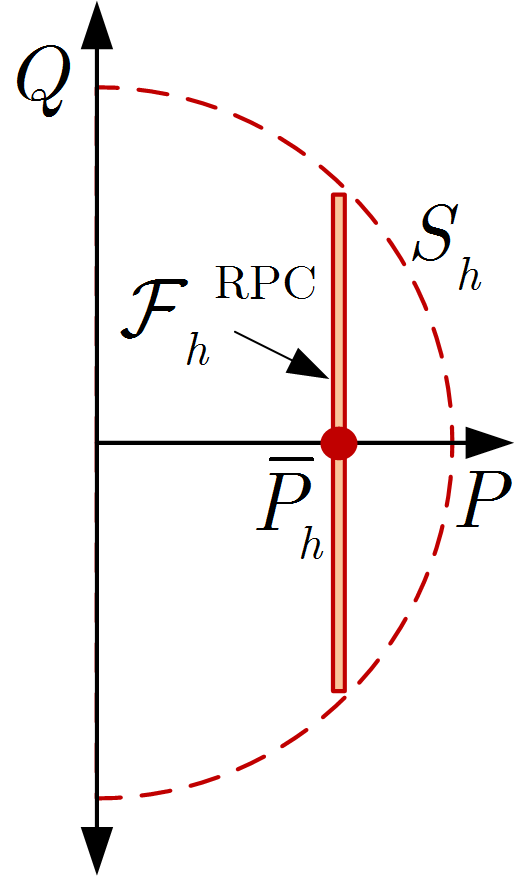}}
  \subfigure[]{\includegraphics[width=0.10\textwidth]{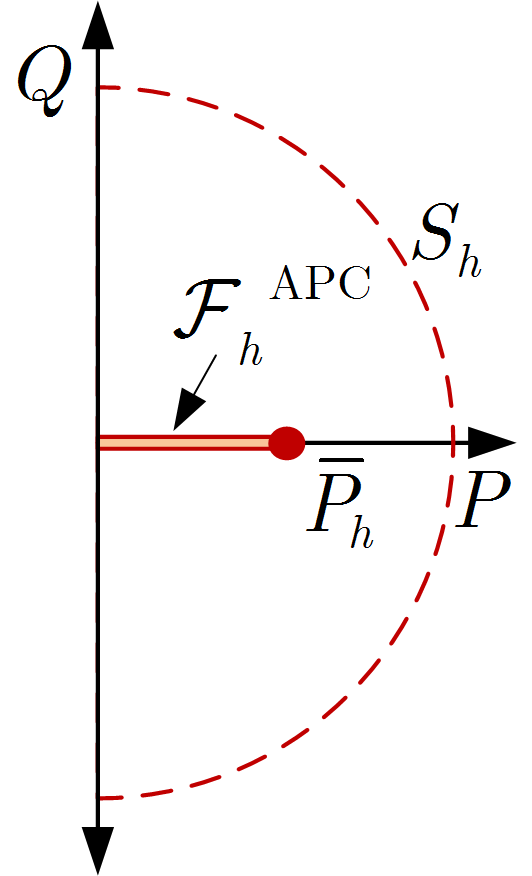}}
  \subfigure[]{\includegraphics[width=0.10\textwidth]{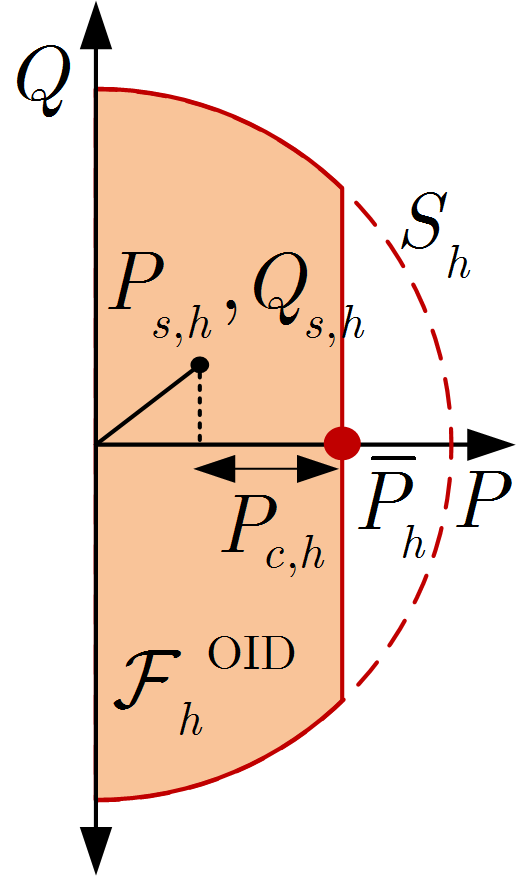}} 
  \subfigure[]{\includegraphics[width=0.10\textwidth]{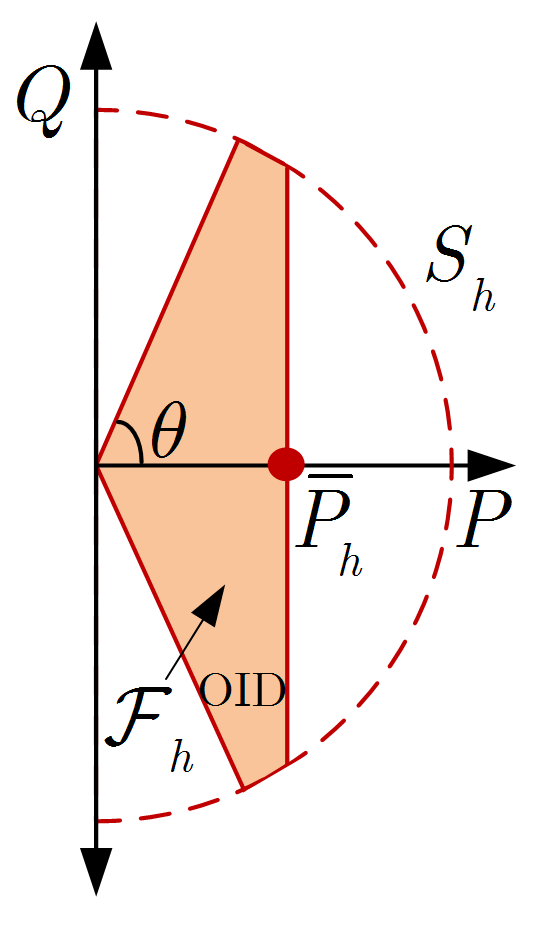}} 
  \caption{Operating regions for the PV inverters are shown by the shaded regions for different strategies: (a) RPC; (b) APC at unity PF; (c) proposed OID strategy; and, (d)
proposed OID strategy with a constraint on the minimum PF.} 
\label{F_PV_inverter_model} \vspace{-.4cm}
\end{figure}

The maximum reactive power that can be generated/consumed by the $h$-th inverter is limited by its rated apparent power, which is denoted here as $S_h$, and the available active power $\bar{P}_h$~\cite{Liu08,Turitsyn11}. Specifically, given $S_h$, the inverter operating space in the $PQ$-plane is given by the semicircle $\{(P,Q): 0 \leq P \leq S_{h}, |Q| \leq \sqrt{S_{h}^2 - P^2}\}$, whose contour is represented by a red dotted line in Fig.~\ref{F_PV_inverter_model}. 

If the $h$-th PV inverter allows reactive power control (RPC), the set of its operating points is given by
\begin{align}
\hspace{-.4cm} \cF_h^\mathrm{RPC} :=  \{(P_{s,h},Q_{s,h}): P_{s,h} = \bar{P}_h,  |Q_{s,h}| \leq \sqrt{S_{h}^2 - \bar{P}_h^2}\}\hspace{-.2cm} 
\label{eq:reactiveCompensation} 
\end{align}
which indicates that the active power output is the available PV power, and the reactive power capability is limited by the inverter rating~\cite{Liu08,Braun10,Turitsyn11}. The set $\cF_h^\mathrm{RPC}$ is represented by the vertical segment in Fig.~\ref{F_PV_inverter_model}(a). When the inverter is not oversized~\cite{Liu08}, no reactive compensation is possible when $\bar{P}_h = S_{h}$; on the other hand, the entire inverter rating can be utilized to supply reactive power when no active power is produced.

Similarly, the set of feasible operating points when an active power curtailment (APC) strategy is implemented at inverter $h$~\cite{DeBrabandere04,Tonkoski11,Tonkoski11Renewable}, is given by
\begin{equation}
\cF_h^\mathrm{APC} := \left\{(P_{s,h},Q_{s,h}):  0 \leq P_{s,h} \leq \bar{P}_h, Q_{s,h} = 0 \right\}
\label{eq:activePowerCurtailment} 
\end{equation}
and it is depicted in Fig.~\ref{F_PV_inverter_model}(b).

In this paper, an OID strategy whereby PV inverters are allowed to adjust \emph{both} active and reactive powers is proposed. Consequently, the set of possible operating points is given by
\begin{align}
\cF_h^\mathrm{OID} := \Big\{(P_{s,h},Q_{s,h}): \,\,&0 \leq P_{s,h} \leq \bar{P}_h, \nonumber \\
& |Q_{s,h}| \leq \sqrt{S_{h}^2 - P_{s,h}^2} \Big\} \, . 
\label{eq:combinedStrategy} 
\end{align}
As shown in Fig.~\ref{F_PV_inverter_model}(c), set $\cF_h^\mathrm{OID}$ is significantly larger then $\cF_h^\mathrm{RPC}$ and $\cF_h^\mathrm{APC}$, thus offering improved flexibility when optimizing the performance of the residential distribution system. Clearly, the control strategy~\eqref{eq:combinedStrategy} subsumes RPC and APC.
A variant of the OID strategy is when the PV inverters are required to operate at a sufficiently high PF to limit the circulation of reactive power throughout the network (see e.g.,~\cite{Braun10}). Thus, if the minimum allowed PF is prescribed to be $\cos\theta$, then the amount of reactive power injected or absorbed is constrained by the two bounds: $|Q_{s,h}| \leq \sqrt{S_{h}^2 - P_{s,h}^2}$  
and $\cos\theta \leq P_{s,h}/ \sqrt{P_{s,h}^2 + Q_{s,h}^2}$. 

The resultant set of feasible inverter operational points is illustrated in Fig.~\ref{F_PV_inverter_model}(d). 

It is worth mentioning that a key implementation task relates to formulating active- and reactive-power setpoints for PV systems composed of multiple inverters, e.g., in microinverter systems. Commercially available microinverter systems usually include a power management unit that communicates via power-line communication (see, e.g., \cite{SolarBridge,Enphase}). Consequently, the optimal real- and reactive-power setpoints could be relayed from the power management units to the different inverters via the power line communication network.

\section{Optimal Inverter Dispatch Problem Formulation}
\label{sec:Dispatch}

Given the operational voltage limits in the residential network, the demanded loads $\bp_{\ell}, \bq_{\ell}$, and the available PV power at the different houses $\bar{\bp}_h$, the objective is to \emph{optimally dispatch} the PV inverters, which amounts to  finding the real and reactive power operating points $\bp_s, \bq_s$  that optimize the operation of the network according to a well defined criterion. To this end, pertinent cost functions are discussed next, while the OID problem is formulated in the subsequent Section~\ref{sec:OptimizationProblem}.

\subsection{Cost Function Formulation}
\label{sec:cost}
The cost function to be minimized is defined as
\begin{equation}
\kappa(\bv, \bi, \bp_c, \bq_s) := c_\rho \rho(\bv,\bi) +  c_\phi \phi(\bp_c) + c_\nu \nu(\bv)  \label{eq:CostFunction}
\end{equation}
where function $\rho(\bv,\bi)$ captures real power losses in the network; $\phi(\bp_c)$ models possible costs of curtailing active power; and, $\nu(\bv)$ promotes a flat voltage profile. Finally, $c_\rho > 0,\,c_\phi \geq 0 ,\,c_\nu \geq 0$, denote weighting coefficients. The three components in~\eqref{eq:CostFunction} are substantiated next.  

\subsubsection{Power losses in the network}
Active power losses in the distribution network are given by~\cite{GSO08} 
\begin{equation}
\rho(\bv,\bi):=\sum_{(m,n) \in \cE} \Re\{V_m I_{mn}^*\} - \Re\{V_n I_{mn}^*\}
\end{equation}
where $I_{mn} \in \mathbb{C}$ denotes the current flowing on line $(m,n)$. 

\subsubsection{Cost associated with active power set points}
The cost incurred by the utility when inverter $h$ is required to operate on a set point different than $(\bar{P}_h,0)$, is modeled without loss of generically by the following quadratic function
\begin{equation}
\label{eq:phiP}
\phi(\bp_{c}) :=  \sum_{h \in \cH} a_h P_{c,h}^2 + b_h P_{c,h} 
\end{equation}
where the choice of the coefficients is based on specific utility-customer prearrangements~\cite{Braun10}. For example, if economic indicators are of interest, the coefficients $a_h$ and $b_h$ could be based on the price at which electricity generated by the PV systems would be sold back to the utility. Similarly, the amount of curtailed power can be minimized by choosing $b_h = 1$ and $a_h = 0, \forall \,\, h \in \cH$.  

\subsubsection{Voltage deviations from average}
Voltage deviations throughout the network can be minimized by defining  
\begin{equation}
\nu(\bv) := \sqrt{\sum_{n \in \cN} \left(|V_n|^2 - \frac{1}{N+1} \sum_{i \in \cN } |V_i|^2 \right)^2}.
\end{equation}
Function $\nu(\bv)$ encourages flat voltage profiles, since it captures the distance of the vector collecting $\{|V_n|^2\}_{n \in \cN}$ from the \emph{average vector} $(1/(N+1) \sum_{i \in \cN} |V_i|^2) \mathbf{1}_{N+1}$. It is important to note that voltage limits will be enforced in the OID problem even if deviations from average are not penalized (that is, when $c_\nu = 0$). Nevertheless, as will be demonstrated in the case studies, including this term provides added flexibility in
reducing deviations of the network voltages from the average.

\subsection{OID Problem} 
\label{sec:OptimizationProblem}
Given the cost function in~\eqref{eq:CostFunction}, the OID problem is formulated as follows:
\begin{subequations} 
\label{Pmg}
\begin{align} 
 \mathrm{(P1)}  \hspace{1.8cm} & \hspace{-1.7cm} \min_{\substack{\{V_n, I_n\} \\ \{P_{c,n},Q_{s,n}, x_n\}}} \,\, \kappa(\bv, \bi, \bp_c, \bq_s) \label{mg-cost} \\
\mathrm{subject\,to} \,\, & \eqref{eq:OhmLaw}, \mathrm{and}  \nonumber  \\ 
\Re\{V_h I_h^*\} & = - P_{\ell,h} + \bar{P}_{h} - P_{c,h}  \hspace{.65cm} \forall \, h \in \cH \label{mg-balance-P} \\
\Im\{V_h I_h^*\} & = - Q_{\ell,h} + Q_{s,h}   \hspace{1.4cm}  \forall \, h \in \cH \label{mg-balance-Q} \\ 
V_u I_u^* & = 0   \hspace{3.3cm} \forall \, u \in \cU  \label{mg-balance-pole} \\ 
V^{\mathrm{min}} & \leq |V_n| \leq V^{\mathrm{max}}  \hspace{1.58cm}  \forall \, n \in \cN  \hspace{-.2cm} \label{mg-Vlimits} \\
0 & \leq P_{c,h}  \leq x_h \bar{P}_{h} \hspace{1.6cm} \forall \, h \in \cH \label{mg-PVp} \\
Q_{s,h}^2 & \leq x_h (S_{h}^2 - (\bar{P}_{h} - P_{c,h})^2) \hspace{.15cm} \forall \, h \in \cH \label{mg-PVq} \\
Q_{s,h} & \leq x_h \tan \theta (\bar{P}_h - P_{c,h}) \hspace{0.65cm} \forall \, h \in \cH \label{mg-pf1} \\
- Q_{s,h} & \leq x_h \tan \theta (\bar{P}_h - P_{c,h}) \hspace{0.65cm} \forall \, h \in \cH \label{mg-pf2} \\
\sum_{h \in \cH} x_h & \leq K \, , \quad \{x_h\} \in \{0,1\}^{|\cH|}  \label{mg-sn}
\end{align}
\end{subequations}
where $x_h$ is a binary variable indicating whether PV-inverter $h$ is controlled ($x_h = 1$) or not ($x_h = 0$), and it is assumed that $K < |\cH|$ PV-inverters are to be controlled. Clearly, when $x_h = 0$, it follows from~\eqref{mg-PVp} and~\eqref{mg-PVq} that $P_{c,h} = 0$ and $Q_{s,h} = 0$, thus implying that PV-inverter $h$ is delivering the maximum available PV power, $\bar{P}_{h}$, at unity power factor. Notice that, when $x_h = 1$, the set of admissible setpoints $\cF_h^\mathrm{OID}$ for inverter $h$ is specified by constraints~\eqref{mg-PVp}--\eqref{mg-pf2} [cf. Figs.~\ref{F_PV_inverter_model}(c) and~\ref{F_PV_inverter_model}(d)]. Finally, the constraint on $V_0$ is left implicit. 

The constraint in~\eqref{mg-sn} offers the possibility of choosing (or minimizing) the number of controlled inverters in order to contain the net operational cost of the residential feeder. Additionally, alternating between different subsets of $K < |\cH|$ facilitates a more equitable treatment of the inverters (see e.g.,~\cite{Turitsyn11}). Particularly important is constraint~\eqref{mg-Vlimits}, whose goal is to ensure the nodal voltages lie between $V^{\mathrm{min}}$ and $V^{\mathrm{max}}$ (even during intervals with peak solar power generation and low demand). Lastly,~\eqref{mg-pf1}-\eqref{mg-pf2}  follow from the PF constraint explained in Section~\ref{sec:Inverter}. 

For a fixed assignment of the binary variables $\{x_h\}_{h \in \cH}$, $(\mathrm{P}1)$ boils down to an instance of the OPF problem. Indeed, similar to various OPF renditions, $(\mathrm{P}1)$ is a nonlinear \emph{nonconvex} problem because of the balance equations~\eqref{mg-balance-P}-\eqref{mg-balance-pole}, the bilinear terms in~\eqref{eq:CostFunction}, and the lower bound in~\eqref{mg-Vlimits}. Furthermore, the presence of the binary variables $\{x_h\}_{h \in \cH}$ renders $(\mathrm{P}1)$ an \emph{NP-hard} MINP, and finding its global optimum requires solving a number of subproblems that increases exponentially ($2^{|\cH|}$) in the number of PV inverters. In principle, off-the-shelf MINP solvers can be employed to find a solution to $(\mathrm{P}1)$; see e.g.,~\cite{PaudyalyISGT}. However, since these solvers are computationally burdensome, they are not suitable for real-time feeder optimization~\cite{PaudyalyISGT}. Further, they do not generally guarantee global optimality of their solutions, which translates here to higher power losses and operational costs. 

\section{Convex Reformulation}
\label{sec:ConvexDispatch}
In this section, a computationally affordable
convex reformulation of $(\mathrm{P}1)$ is derived by leveraging
sparsity-promoting regularizations to drop the binary selection
variables (Section~\ref{sec:sparsity}), and utilizing semidefinite relaxation
techniques to address the nonconvexity due to bilinear terms and voltage lower
bounds (Section~\ref{sec:sdp}).

\subsection{Sparsity-leveraging OID}
\label{sec:sparsity}

If PV-inverter $h$ is not selected, then it provides the maximum active power $\bar{P}_h$ at unity PF; consequently, one clearly has that $P_{c,h} = Q_{s,h} = 0$. On the other hand, $P_{c,h}$ and $Q_{s,h}$ may be different than zero with OID [cf. Fig.~\ref{F_PV_inverter_model}]. 
Supposing that $K < |\cH|$, it follows that the $2 |\cH|\times 1$ real-valued vector $[\bp_c^\cT, \bq_s^\cT]^\cT$ is \emph{group sparse}~\cite{YuLi06}; meaning that, the $2 \times 1$ sub-vectors $[P_{c,h}, Q_{s,h}]^\cT$ (i.e., the ``groups'' of variables~\cite{YuLi06}) corresponding to inverters operating at point $(\bar{P}_h,0)$ contain all zeroes, whereas $[P_{c,h}, Q_{s,h}]^\cT \neq \mathbf{0}_{2 \times 1}$ for the inverters operated under OID. For example, consider the system in Fig.~\ref{F_LV_network}, and suppose that only PV-inverters in $\mathrm{H}_{10}, \mathrm{H}_{12}$ are to be controlled. Then, it follows that $[\bp_c^\cT, \bq_s^\cT]^\cT = [\mathbf{0}_{18 \times 1}^{\cT}, P_{c,10}, Q_{s,10}, \mathbf{0}_{2 \times 1}^{\cT}, P_{c,12}, Q_{s,12}]^\cT$, with $[P_{c,10}, Q_{s,10}]^\cT \neq \mathbf{0}_{2 \times 1}$ and $[P_{c,12}, Q_{s,12}]^\cT \neq \mathbf{0}_{2 \times 1}$.

One major implication of this group sparsity attribute of the real-valued vector $[\bp_c^\cT, \bq_s^\cT]^\cT$, is that one can discard the binary variables $\{x_{h}\}_{h \in \cH}$, and effect PV inverter selection by leveraging sparsity-promoting regularization techniques. Among possible candidates, the following group-Lasso-type function is well suited for the problem at hand~\cite{YuLi06,Wiesel11}:
\begin{align}
\gamma(\bp_c, \bq_s) :=  \lambda \sum_{h \in \cH} \,  \sqrt{ P_{c,h}^2 + Q_{s,h}^2}  \label{Glasso_powers}  
\end{align}
where $\lambda \geq 0$ is a tuning parameter. Thus, using~\eqref{Glasso_powers}, the OID problem in $(\mathrm{P}1)$ can be relaxed to:
\begin{subequations} 
\label{Pmg2}
\begin{align} 
 \mathrm{(P2)}  \hspace{1.8cm} & \hspace{-1.7cm} \min_{\substack{\{V_n, I_n\} \\ \{P_{c,n},Q_{s,n}\}}} \,\, \kappa(\bv, \bi, \bp_c, \bq_s) + \gamma(\bp_c, \bq_s) \label{mg-cost2} \\
\mathrm{subject\,to} \,\, & \eqref{eq:OhmLaw},~\eqref{mg-balance-P}-\eqref{mg-Vlimits},  \mathrm{and} \, \forall \, h \in \cH \nonumber  \\ 
0 & \leq P_{c,h}  \leq  \bar{P}_{h} \label{mg-PVp2} \\
Q_{s,h}^2 & \leq  S_{h}^2 - (\bar{P}_{h} - P_{c,h})^2  \hspace{.2cm} \label{mg-PVq2} \\
Q_{s,h} & \leq \tan \theta (\bar{P}_h - P_{c,h}) \hspace{1.2cm} \label{mg-pf1p2} \\
- Q_{s,h} & \leq \tan \theta (\bar{P}_h - P_{c,h}). \hspace{1.2cm} \label{mg-pf2p2}
\end{align}
\end{subequations}
Notice the absence of the binary selection variables $\{x_h\}_{h \in \cH}$ in this formulation. The role of $\lambda$ is to control the number of sub-vectors 
$[P_{c,h}, Q_{s,h}]^\cT$ that are set to $ \mathbf{0}_{2 \times 1}$. Particularly, for $\lambda = 0$, all the inverters may curtail active power and provide reactive power compensation (that is, $K = |\cH|$); and, as $\lambda$ is increased, the number of inverters participating in OID decreases~\cite{YuLi06}. Appendix~\ref{sec:appB} elaborates further on the PV-inverter selection capability offered by the regularization function $\gamma(\bp_c, \bq_s)$. A generalization of~\eqref{Glasso_powers} is represented by the weighted counterpart $\gamma_w(\bp_c, \bq_s) :=  \sum_{h \in \cH} \lambda_h   \sqrt{ P_{c,h}^2 + Q_{s,h}^2}$, where $\{\lambda_{h} \geq 0\}_{h \in \cH}$ substantiate possible preferences to use (low value of $\lambda_{h}$) or not (high value of $\lambda_{h}$) specific inverters. With this approach, a more equitable treatment of PV inverters can be facilitated by discouraging the use of inverters whose active and reactive powers are adjusted more frequently.

\subsection{SDP reformulation}
\label{sec:sdp}
Key to deriving an SDP reformulation of $(\mathrm{P}2)$ is to express the active and reactive powers injected per node, and voltage magnitudes as \emph{linear} functions of the outer-product matrix $\bV := \bv \bv^\cH$. 
First, let us consider the constraints in $\mathrm{(P2)}$. Let $\{\be_n\}_{n = 1}^{N+1}$ denote the canonical basis of $\mathbb{R}^{N+1}$, and define the admittance-related matrix $\bY_n := \be_n \be_n^\cT \bY$ per node $n \in \cN$. Furthermore, define the Hermitian matrices $\bZ_n := \frac{1}{2}(\bY_n + \bY_n^\cH)$, $\bar{\bZ}_n := \frac{j}{2}(\bY_n - \bY_n^\cH)$, and $\bM_n := \be_n \be_n^\cT$. Then, the nonconvex constraints~\eqref{mg-balance-P}-\eqref{mg-balance-pole} and~\eqref{mg-Vlimits} can be equivalently re-written as (see~\cite{LavaeiLow,Tse12,Dallanese-TSG13} for further details):
\begin{subequations} 
\label{SDPconst}
\begin{align} 
\trace(\bZ_h \bV) & =  - P_{\ell,h} + \bar{P}_{h} - P_{c,h},  & \forall \, h \in \cH \label{SDPconstP} \\
\trace(\bar\bZ_h \bV) & =  - Q_{\ell,h} + Q_{s,h},  & \forall \, h \in \cH \label{SDPconstQ} \\
\trace(\bZ_n \bV) & =  0 , \,\, \trace(\bar{\bZ}_n \bV) =  0,  & \forall \, n \in \cU \label{SDPconstPole} \\
(V^{\mathrm{min}})^2 & \leq \trace(\bM_n \bV) \leq (V^{\mathrm{max}} )^2,  &  \forall \, n \in \cN \label{SDPconstV}
\end{align}
\end{subequations}
which are all linear in the matrix variable $\bV$.

Now, let us consider the cost function in~\eqref{eq:CostFunction}. Define the projector matrix $\bPi := \bI_{N+1} - \frac{1}{N+1}\mathbf{1}_{N+1 \times 1} \mathbf{1}_{N+1 \times 1} ^\cT$, and let $\diag(\bV)$ denote the operator that returns a vector collecting the diagonal elements of matrix $\bV$. Then, the function $\nu(\bv)$ in~\eqref{eq:CostFunction} can be equivalently expressed in terms of $\bV$ as $\nu(\bV) := \| \bPi \diag(\bV) \|_2$.  Next, for each distribution line $(m,n) \in \cE$, define the symmetric matrix $\bL_{mn} := \Re\{y_{mn}\} (\be_m - \be_n) (\be_m - \be_n)^\cT$, and notice that the active power loss on the line $(m,n) \in \cE$ can be re-expressed as a linear function of the matrix $\bV$ as $ \Re\{V_m I_{mn}^*\} - \Re\{V_n I_{mn}^*\} = \trace(\bL_{mn} \bV)$. Thus, using matrices $\bPi$ and $\bL := \sum_{(m,n) \in \cE} \bL_{mn}$, the cost function in~\eqref{eq:CostFunction} becomes
\begin{align} 
\hspace{-.2cm} \kappa(\bV, \bp_c, \bq_s) &= c_\rho \trace(\bL \bV) + c_\phi \phi(\bp_c) + c_\nu \| \bPi \diag(\bV) \|_2 .  \hspace{-.1cm} 
\label{eq:CostFunction2}
\end{align}
Notice that for $c_\nu \geq 0$, the surrogate cost $\kappa(\bV, \bp_c, \bq_s)$ is convex in the matrix variable $\bV$. 

Aiming for an SDP formulation of $(\mathrm{P}2)$, Schur's complement can be leveraged to convert the non-linear summands in~\eqref{Glasso_powers} and~\eqref{eq:CostFunction2} to a linear cost over auxiliary optimization variables~\cite{Vandenberghe96} (see Appendix~\ref{sec:appA} for details). Specifically, with $\{y_h, w_h, h = 1,\ldots, |\cH|\}_{h \in \cH}$ and $\upsilon$ serving as auxiliary variables, $(\mathrm{P}2)$ can be equivalently reformulated as
\begin{subequations} 
\label{Pmg3}
\begin{align} 
 &\hspace{-.1cm} \mathrm{(P3)}  \min_{\substack{\{P_{c,n},Q_{s,n}\} \\ \{y_h, w_h \geq 0 \} \\ \bV , \upsilon } } \hspace{-0.2cm} c_\rho \trace(\bL \bV) + c_\phi  \hspace{-.1cm} \sum_{h \in \cH} y_h +  c_\nu \upsilon +  \lambda \hspace{-.1cm} \sum_{h \in \cH}w_h \nonumber \\
& \mathrm{subject\,to} \,\,\eqref{mg-PVp2},~\eqref{mg-pf1p2}-\eqref{mg-pf2p2},~\eqref{SDPconst},  \mathrm{and} \, \forall \, h \in \cH  \nonumber  \\ 
& \left[ 
\begin{array}{cc}
\upsilon \bI_{N+1} & \bPi \diag(\bV) \\
\diag(\bV)^{\cT} \bPi^\cT & \upsilon 
\end{array}
\right] \succeq \mathbf{0}  \label{SDPvoltage} \\
& \left[ \begin{array}{cc}
b_h P_{c,h} - y_h & \sqrt{a_h} P_{c,h} \\
\sqrt{a_h} P_{c,h} & - 1 
\end{array}
\right] \preceq \mathbf{0} \hspace{1.5cm} 
\label{SDPConstCostActive} \\
& \left[ \begin{array}{ccc}
w_h & 0 & P_{c,h} \\
0 & w_h & Q_{s,h} \\
P_{c,h} & Q_{s,h} & w_h
\end{array}
\right] \succeq \mathbf{0} \hspace{2.1cm} 
\label{SDPnorm2} \\
& \left[ 
\begin{array}{ccc}
-S_h^2 & Q_{s,h} & \bar{P}_{h} - P_{c,h} \\
Q_{s,h} & -1 & 0 \\
\bar{P}_{h} - P_{c,h}  & 0 & -1
\end{array}
\right] \preceq \mathbf{0} \hspace{.5cm}
\label{SDPConstReactive} \\
& \,\, \bV \succeq \mathbf{0}  \label{SDPpositive} \\
&  \,\, \rank(\bV) = 1 \,.   \label{SDPRank}
\end{align}
\end{subequations}
The positive semi-definiteness and rank constraints~\eqref{SDPpositive}--\eqref{SDPRank} jointly ensure that for any feasible $\bV \in \mathbb{C}^{N+1 \times N+1}$, there always exists a vector of voltages $\bv$ such that $\bV = \bv \bv^\cH$, and that formulations $(\mathrm{P}2)$ and $(\mathrm{P}3)$ solve an equivalent OID problem. For the optima of these two problems, it holds that $\bV^\mathrm{opt} = \bv^\mathrm{opt} (\bv^\mathrm{opt})^\cT$; $\upsilon^\mathrm{opt} =  \| \bPi \diag(\bV^\mathrm{opt}) \|_2$; $\sum_{h \in \cH} y_h^\mathrm{opt} + z_h^\mathrm{opt} = \phi( \bp_{c}^\mathrm{opt})$;  $\trace(\bL \bV^\mathrm{opt}) = \rho(\bv^\mathrm{opt},\bi^\mathrm{opt})$; and $\lambda \sum_{h \in \cH} w_h^\mathrm{opt} = \gamma(\bp_c^\mathrm{opt}, \bq_s^\mathrm{opt})$, where superscript $(\cdot)^\mathrm{opt}$ is used to denote the globally optimal solution. 

Unfortunately, $(\mathrm{P}3)$ is still nonconvex due to the rank-$1$ constraint in~\eqref{SDPRank}. Nevertheless, in the spirit of the semidefinite relaxation technique~\cite{luospmag10}, one can readily obtain the following (relaxed) SDP reformulation of the OID problem by dropping the rank constraint:
\vspace{-.1cm}
\begin{align} 
& \hspace{-.3cm} \mathrm{(P4)}  \hspace{-.4cm} \min_{\substack{\{P_{c,n},Q_{s,n}\} \\ \{y_h, w_h \geq 0 \} \\ \bV , \upsilon } } \hspace{-0.4cm} c_\rho \trace(\bL \bV) + c_\phi  \hspace{-.1cm} \sum_{h \in \cH}(y_h + z_h) +  c_\nu \upsilon +  \lambda \hspace{-.1cm} \sum_{h \in \cH}w_h  \nonumber \\
& \mathrm{subject\,to} \,\,  \eqref{mg-PVp2},~\eqref{mg-pf1p2}-\eqref{mg-pf2p2} \nonumber \\
& \hspace{1.4cm}  ~\eqref{SDPconst},~\eqref{SDPvoltage}-\eqref{SDPConstReactive},~\mathrm{and} \,\, \bV \succeq \mathbf{0}.  
\end{align}
If the optimal solution $\bV^\mathrm{opt}$ of $(\mathrm{P}4)$ has rank 1, then it is also a \emph{globally} optimal solution for the nonconvex problems $(\mathrm{P}3)$ and $(\mathrm{P}2)$~\cite{luospmag10}. Further, given $\bV^\mathrm{opt}$, the vectors of voltages and injected currents can be computed as $\bv^\mathrm{opt} := \sqrt{\lambda_1} \bu_1$ and $\bi^\mathrm{opt} = \bY \bv^\mathrm{opt} $, respectively, where $\lambda_1 \in \mathbb{R}^+$ is the unique non-zero eigenvalue of $\bV^\mathrm{opt}$, and $\bu_1$ the corresponding eigenvector.   

A caveat of the semidefinite relaxation technique is that matrix $\bV$ could have rank greater than $1$; in this case, the feasible rank-1 approximation of $\bV$ obtainable via rank reduction techniques turns out to be generally suboptimal~\cite{luospmag10}. Fortunately though, when the considered radial residential feeder is single-phase, and the cost~\eqref{eq:CostFunction2} is strictly increasing in the line power flows, Theorem~1 of~\cite{Tse12} can be conveniently adapted to the problem at hand to show that a rank-$1$ solution of $(\mathrm{P}4)$ is always attainable provided that $(\mathrm{P}2)$ is feasible, and the controlled inverters are sufficiently oversized. Incidentally, allowing inverters to operate in the region $\cF_h^\mathrm{OID}$ facilitates obtaining a rank-1 solution, since higher reactive powers can be absorbed (a sufficient condition in~\cite[Thm.~1]{Tse12} requires nodes to be able to absorb a sufficient amount of reactive power). 
When the house load is not balanced, the OPF formulation for unbalanced distribution systems in~\cite{Dallanese-TSG13} can be applied to account for discrepancies in nodal voltages and conductor currents.

\vspace{1mm}

\noindent \emph{Remark~1}. The number of controlled PV inverters depends on the parameter $\lambda$.
Provided $(\mathrm{P}2)$ is feasible, there are (at least) three viable ways to select $\lambda$ so that $K$ inverters are selected: i) First, since both solar irradiation and load typically manifest daily/seasonal patterns, $\lambda$ can be selected based on historical OID results. ii) In the spirit of the so-called cross-validation~\cite{YuLi06}, OID problems featuring different values of $\lambda$ can be solved in parallel. iii) Finally, notice that the optimized values of $\bp_c, \bq_s$  are inherently related to the dual variables associated with constraints~\eqref{SDPconstP}--\eqref{SDPconstQ}~\cite{YuLi06,Wiesel11}; then, when a primal-dual algorithm is employed to solve $(\mathrm{P}4)$, $\lambda$ can be conveniently adjusted during the iterations of this scheme.   
\vspace{1mm}

\emph{Remark~2}. A second-order cone program (SOCP) relaxation of the OPF problem for radial systems was recently proposed in~\cite{Farivar13}. However, sufficient conditions for global optimality of the obtained solution require discarding the upper bounds on the real and reactive power injections, which may not be a realistic setup for residential distribution systems. In contrast, conditions derived in~\cite{Tse12} for global optimality of the SDP relaxation are typically satisfied in practice. Furthermore, the virtues of the SDP relaxation in the case of an unbalanced system operation were demonstrated in~\cite{Dallanese-TSG13}.

\subsection{Optimal RPC and APC strategies}
\label{sec:rpc}

It is immediate to formulate optimal RPC and APC strategies in the proposed framework. Towards this end, consider the following regularization function:
\vspace{-.1cm}
\begin{align}
\gamma^\prime(\bp_c, \bq_s) :=  \gamma(\bp_c, \bq_s) +   \lambda_p \|\bp_c\|_1 + \lambda_q \|\bq_s\|_1
\label{SGlasso_powers}  
\end{align}
with $\lambda_p \geq 0$ and $\lambda_q \geq 0$ serving as tuning parameters. As in typical (constrained) sparse linear regression problems~\cite{Tib96}, the second and third terms in the right hand side of~\eqref{SGlasso_powers} promote entry-wise sparsity in the vectors $\bp_c$ and $\bq_s$, respectively. Thus, by replacing $\gamma(\bp_c, \bq_s)$ with $\gamma^\prime(\bp_c, \bq_s)$  in $(\mathrm{P}2)$-$(\mathrm{P}4)$ the following inverter dispatch strategies can be readily implemented: i) $\lambda > 0, \lambda_p = \lambda_q = 0$: inverters can curtail real power and provide/consume reactive power. ii) $\lambda = 0, \lambda_p > 0, \lambda_q = 0, \bq_s = \mathbf{0}$: APC-only is implemented at a subset of PV systems. iii) $\lambda = 0, \lambda_p = 0, \bp_c = \mathbf{0}, \lambda_q > 0$: some of the inverters are optimally selected to inject/absorb reactive power. iv) $\lambda > 0, \lambda_p > 0, \lambda_q > 0$: mixed strategy whereby all three strategies (RPC, APC, and OID) can be implemented.  

Compared to the local strategies in~\cite{DeBrabandere04,Carvalho08,Cagnano11,Turitsyn11,Bolognani13,Aliprantis13,Tonkoski11}, the resultant RPC and APC tasks are cast here as instances of the OPF problem, thus ensuring system-level optimality guarantees; and, furthermore, through the regularization function~\eqref{SGlasso_powers}, they offer the flexibility of selecting the subset of critical PV inverters that must be controlled in order to fulfill optimization objectives and electrical network constraints. An OPF-type RPC strategy was previously proposed in~\cite{Farivar12}; however, inverters providing reactive support were preselected.

\begin{table}[b]
\caption{Single-phase $\pi$-model line parameters}
\vspace{-.2cm}
\begin{center}
\begin{tabular}{l||c|c|c}
& $R$ [$\Omega/\textrm{km}$]  & $L$ [$\textrm{mH}/\textrm{km}$] & $C$ [$\mu\textrm{F}/\textrm{km}$]  \\
\hline
Drop line  & 0.549 & 0.230 & 0.055 \\
Pole-pole line & 0.270 & 0.240 & 0.072 
\end{tabular}
\end{center}
\label{tab:lines}
\end{table}%

\begin{figure}[b]
\begin{center}
\includegraphics[width=0.5\textwidth]{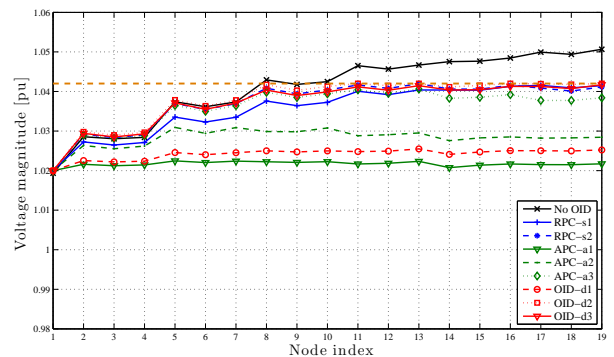}
\vspace{-.5cm}
\caption{Voltage Profile $\{|V_n|\}_{n \in \cN}$ at 12:00 with and without inverter control.}
\label{F_voltage}
\end{center}
\end{figure}

\section{Numerical Case Studies}
\label{sec:CaseStudies}

To implement the OID strategy for residential network optimization, the utility requires: i) the network admittance matrix $\bY$; ii) instantaneous available powers $\{\bar{P}_h\}$; and, iii) the ratings $\{S_h\}$. Once $(\mathrm{P}4)$ is solved, the active- and reactive-power setpoints (i.e., $\bp_s, \bq_s$) are relayed to the inverters. 

Consider the distribution network in Fig.~\ref{F_LV_network}, which is adopted from~\cite{Tonkoski11, Tonkoski12}. The pole-pole distance is set to $50$ m, while the lengths of the drop lines are set to $20$ m. The parameters for the admittance matrix are adopted from~\cite{Tonkoski11, Tonkoski12}, and summarized in Table~\ref{tab:lines}.
Voltages $V^\textrm{min}$ and $V^\textrm{max}$ in $\mathrm{(P1)}$ refer to the minimum and maximum voltage utilization limits for the normal operation of residential systems, and they depend on the specific standard adopted. Since the test system in Fig.~\ref{F_LV_network} is adopted from~\cite{Tonkoski11, Tonkoski12}, these limits are set to 0.917 pu and 1.042 pu, respectively in this case study (see page 11 of the CAN3-C235-83 standard). Notice that typical voltages that dictate inverters to shut down according to~\cite{Standard1547} are higher than $V^\textrm{max}$.

\begin{figure*}[t]
\centering
  \subfigure[]{\includegraphics[width=0.40\textwidth]{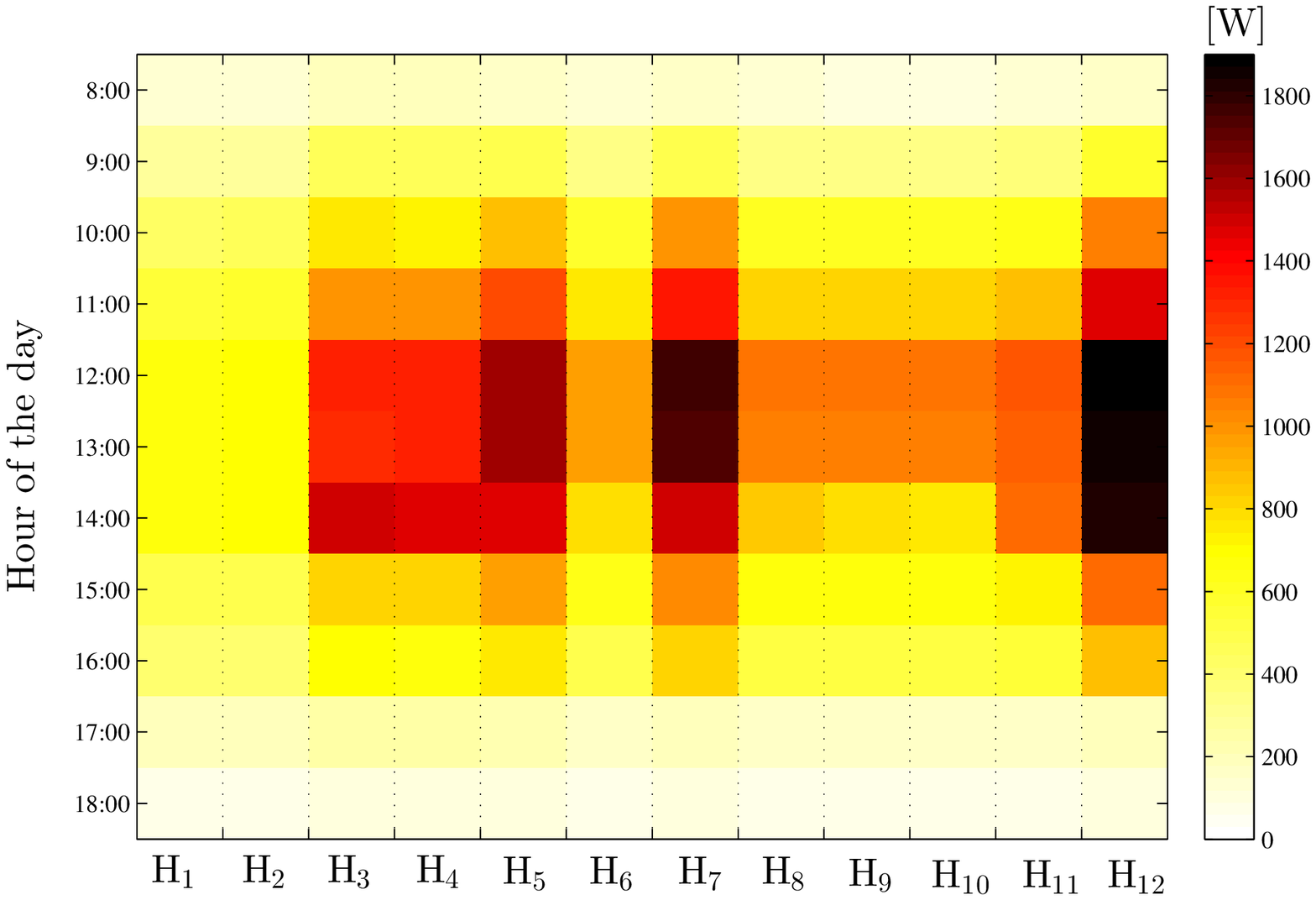}} \vspace{-.2cm}
  \subfigure[]{\includegraphics[width=0.40\textwidth]{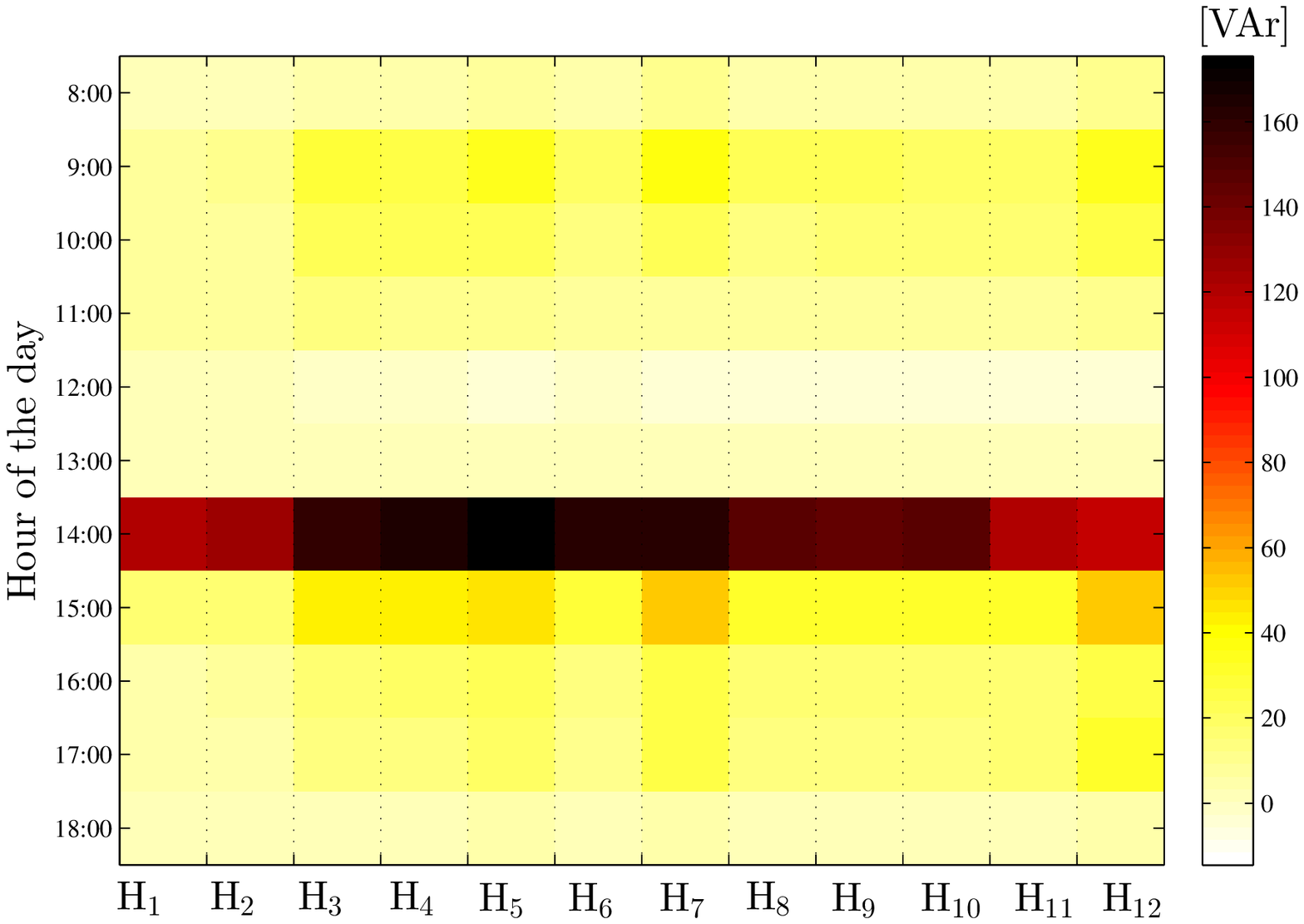}} 
  \subfigure[]{\includegraphics[width=0.40\textwidth]{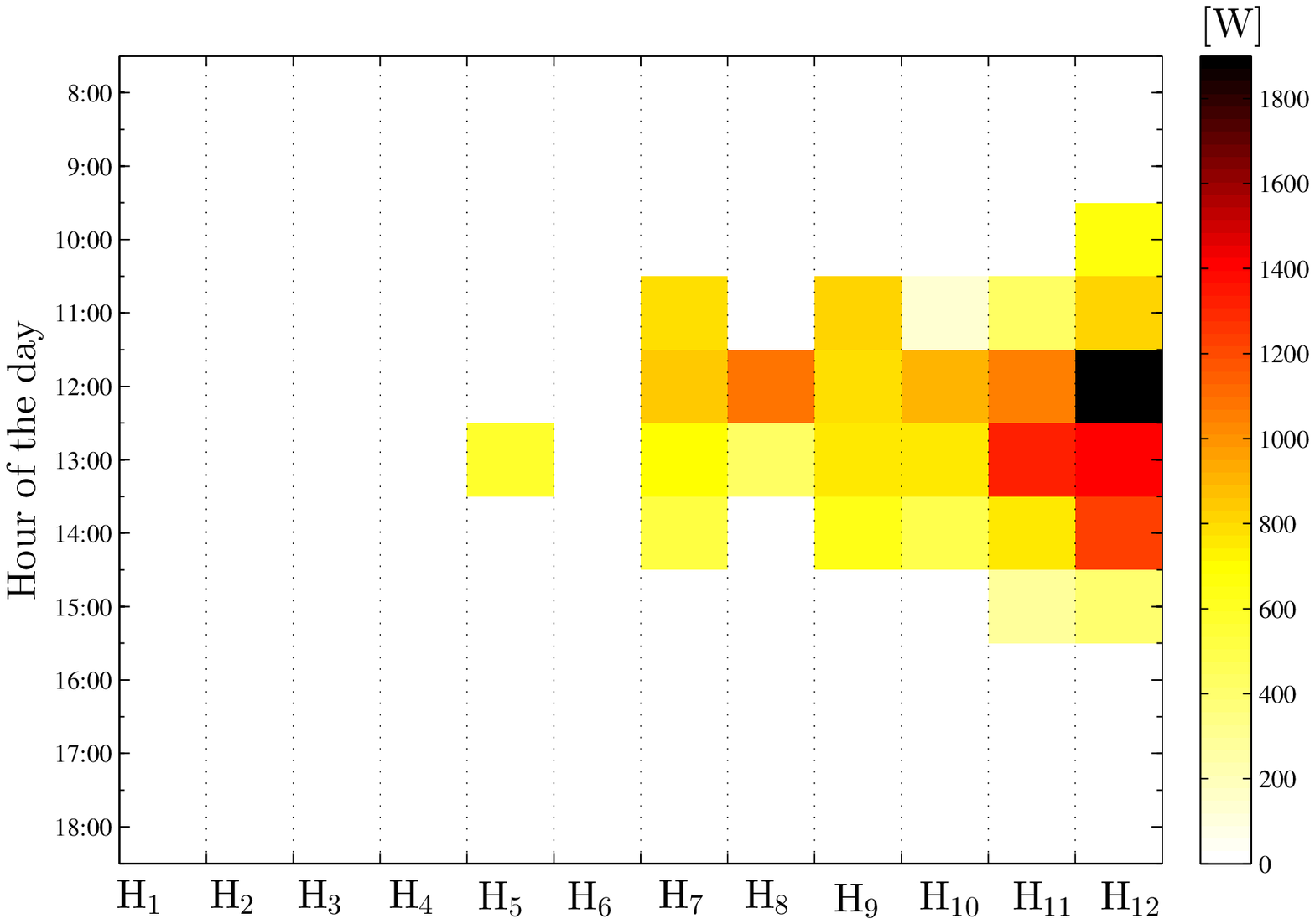}}  
  \subfigure[]{\includegraphics[width=0.40\textwidth]{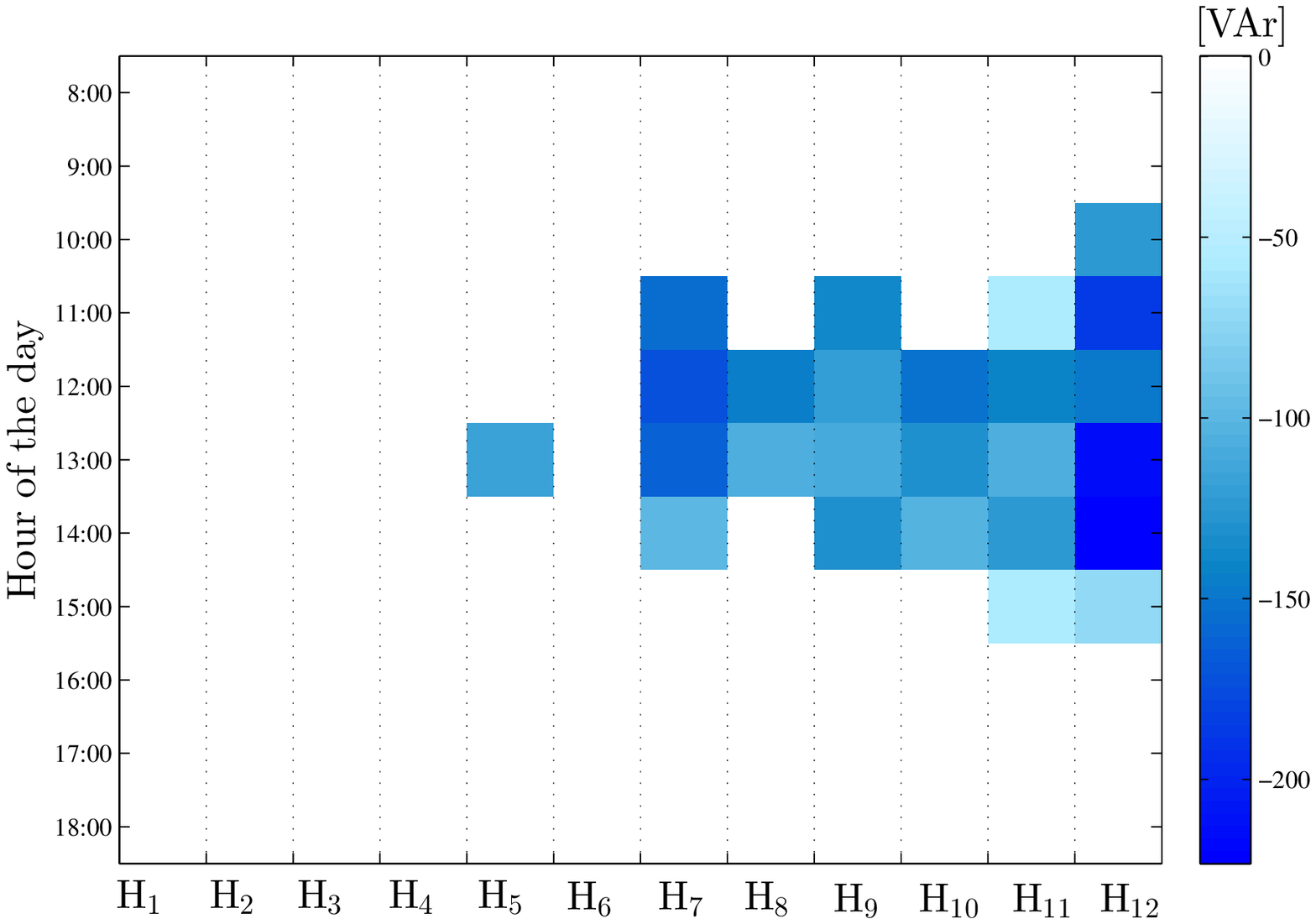}} 
  \caption{Dispatched inverters: Curtailed active power (a) and reactive power (b) for $\lambda = 0$. Curtailed active power (c) and reactive power (d) for $\lambda = 0.8$.} 
\label{F_OID_result} 
\end{figure*}

The optimization package \texttt{CVX}\footnote{[Online] Available: \texttt{http://cvxr.com/cvx/}} along with the interior-point based solver \texttt{SeDuMi}\footnote{[Online] Available: \texttt{http://sedumi.ie.lehigh.edu/}} are employed to solve the OID problem in \texttt{MATLAB}. The average computational time required to solve $(\mathrm{P}4)$ was $0.27$ s on a machine with a Intel Core i7-2600 CPU @ 3.40GHz. In all the conducted numerical tests, the rank of matrix $\bV^\mathrm{opt}$ was always $1$, meaning that the globally optimal solution of (P3) was obtained [cf. Section~\ref{sec:sdp}].

All 12 houses feature fixed roof-top PV systems, with a DC-AC derating coefficient of $0.77$~\cite{King(SANDIA):2007}. The DC ratings of the houses are as follows: $5.52$ kW for houses $\mathrm{H}_1, \mathrm{H}_9, \mathrm{H}_{10}$ (modeled along the lines of~\cite{Cady2012}); $5.70$ kW (which corresponds to the average installed rating for residential systems in 2011~\cite{Sherwood12}) for $\mathrm{H}_2, \mathrm{H}_6, \mathrm{H}_{8}, \mathrm{H}_{11}$; and, $9.00$ kW for the remaining five houses (modeled along the lines of~\cite{Dhople2009}). As suggested in~\cite{Turitsyn11}, it is assumed that the PV inverters are oversized by $10\%$ of the resultant AC rating. 
Further, the minimum PF for the inverters is set to 0.85~\cite{Braun10}. The available active powers $\{\bar{P}_h\}$ during the day are computed using the System Advisor Model (SAM)\footnote{[Online] Available at \texttt{https://sam.nrel.gov/}.} of the National Renewable Energy Laboratory (NREL), based on   
the  typical meteorological year (TMY) data for Minneapolis, MN, during the month of July.
The residential load profile is obtained from the Open Energy Info database,\footnote{``Commercial and Residential Hourly Load Profiles for all TMY3 Locations in the United States,'' accessible from \texttt{http://en.openei.org/datasets/node/961}. Developed at the National Renewable Energy Laboratory and made available under the ODC-BY 1.0 Attribution License.} and the base load experienced in downtown Saint Paul, MN, during the month of July is used for this test case. To generate different load profiles, the base active power profile is perturbed using a Gaussian random variable with zero mean and standard deviation $200$ W and a PF of 0.9 is presumed~\cite{Tonkoski12}.
Finally, the voltage at the secondary of the transformer is set at 1.02 pu, to ensure a minimum voltage magnitude of $0.917$ pu at $\mathrm{H}_{11}-\mathrm{H}_{12}$ when the PV inverters do not generate power. 

With these data, the peak net active power injection (i.e., $\bar{P}_h - P_{\ell,h}$) occurs approximately at solar noon. Standard power-flow computations yield the voltage profile illustrated in Fig.~\ref{F_voltage} with a black solid line. Notice that the voltage magnitude towards the end of the feeder (nodes 11-19) exceeds the upper limit; furthermore, the voltage magnitude at houses $\mathrm{H}_{11}, \mathrm{H}_{12}$ (nodes 17 and 19) is beyond $1.05$ pu, which is the limit usually set for inverter protection~\cite{King(SANDIA):2007}. 

Consider then, implementing the proposed OID strategy, and suppose that the active power losses are to be minimized; that is, 
the weighting coefficients in~\eqref{eq:CostFunction2} are $c_\rho = 1$, $c_\phi = 0$, and $c_\nu = 0$. To emphasize the role of the sparsity-promoting regularization function, Fig.~\ref{F_OID_result} illustrates the active power curtailed and the reactive power provided by each inverter for $\lambda = 0$ (Figs.~\ref{F_OID_result}(a) and (b)) and $\lambda = 0.8$  (Figs.~\ref{F_OID_result}(c) and (d)). For $\lambda = 0$, it is clearly seen that all inverters are controlled; in fact, they all curtail active power from 8:00AM to 6:00PM, and inject reactive power during the entire interval. Interestingly, more active power is curtailed at houses with higher AC ratings ($\mathrm{H}_{3}, \mathrm{H}_{4}, \mathrm{H}_{5}, \mathrm{H}_{7}, \mathrm{H}_{12}$). When $\lambda = 0.8$, only $7$ inverters are controlled, thus corroborating the ability of the regularization function in~\eqref{Glasso_powers} in effecting inverter selection. Four remarks are in order: i) the 7 controlled inverters are the ones located far from the transformer; ii) houses at the end of the feeder curtail more active power and absorb more reactive power than the others (thus matching the findings in~\cite{Tonkoski11, Tonkoski12}); iii) all the selected inverters absorb reactive power (whereas they inject power in Fig.~\ref{F_OID_result}(b)); and iv) it is noted that for $\lambda > 0.8$, no less inverters are selected, meaning that at least $7$ inverters must be controlled in order to effect voltage regulation in this particular case study  [cf.~\eqref{SDPconstV}]. The resultant voltage profiles are illustrated in Fig.~\ref{F_voltage}, where ``OID-d1'' refers to the case $\lambda = 0$, and ``OID-d2'' to the case $\lambda = 0.8$. Clearly, while the voltage limits are enforced
in both cases, a flatter voltage profile is obtained in the second case.

\begin{figure}[t]
\centering
  \subfigure[]{\includegraphics[width=0.42\textwidth]{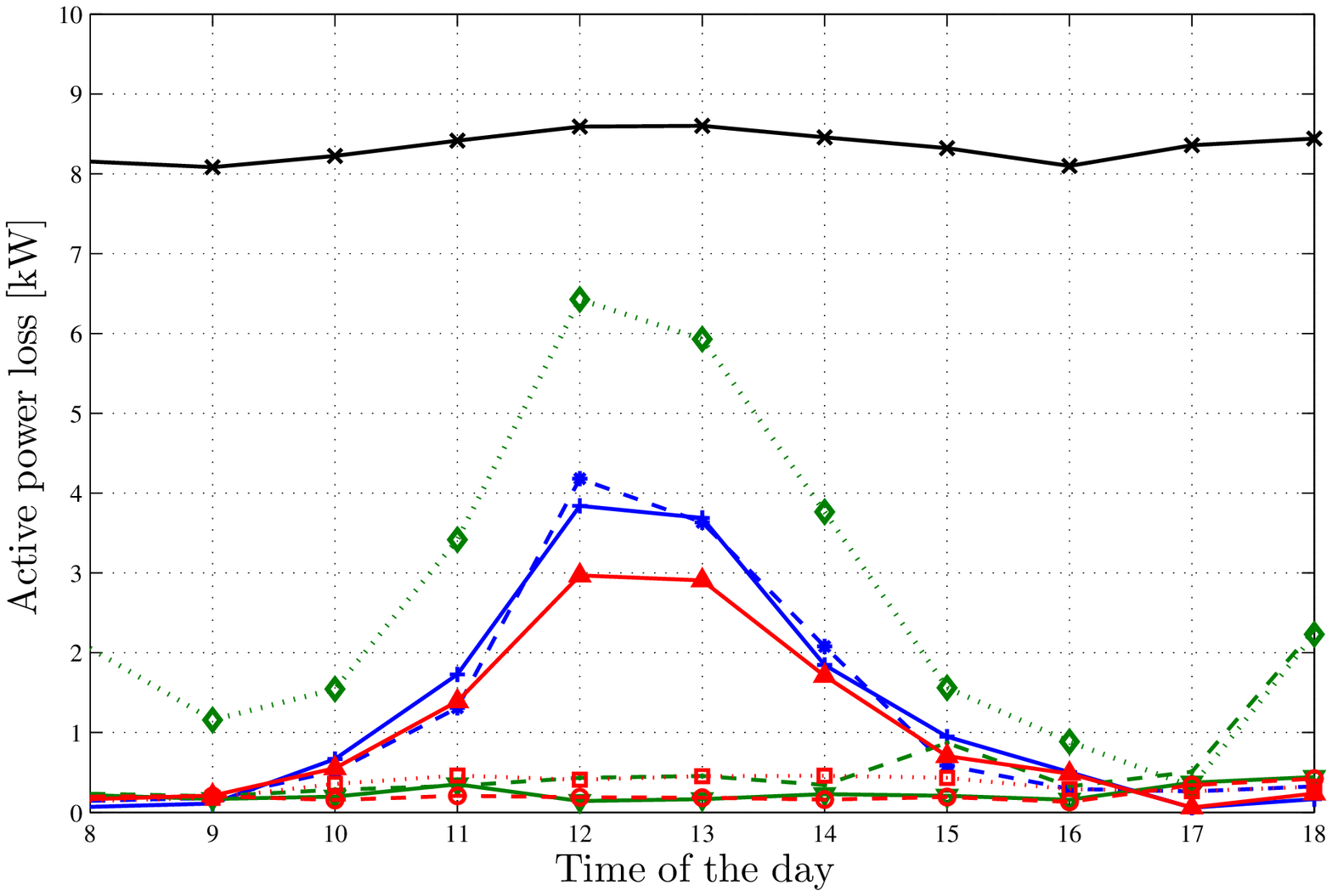} } 
  \subfigure[]{\includegraphics[width=0.42\textwidth]{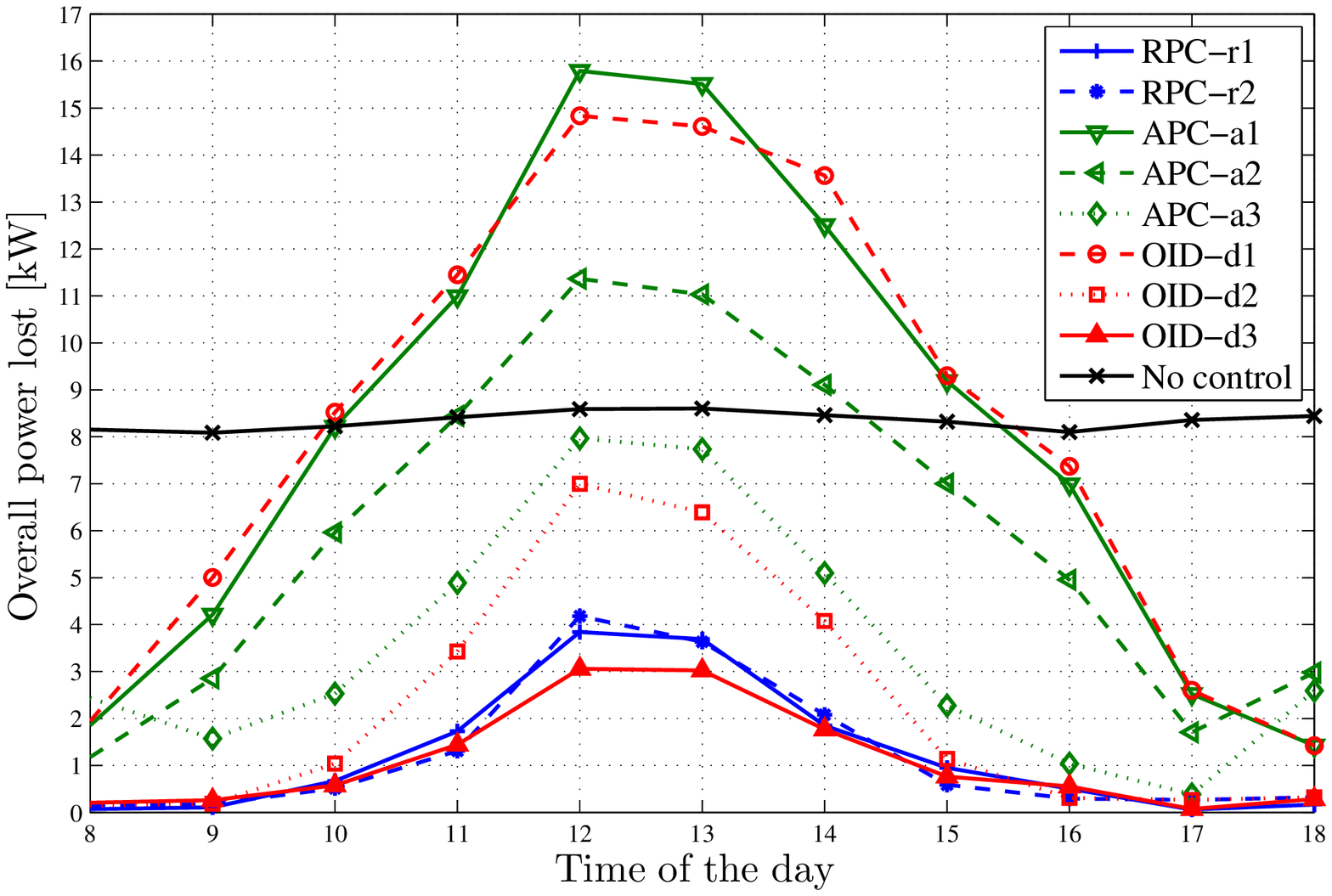}} 
  \caption{(a) Active power loss [kW] on the distribution lines. (b) Sum of the power lost on the distribution lines and the one curtailed by the inverters.} 
\label{F_power_loss} 
\end{figure}

Next, the OID strategy is compared with: i)  RPC without inverter selection (``RPC-r1''); ii) RPC with inverter selection (``RPC-r2''), where $\lambda_q$ is such that only 7 inverters provide reactive powers; iii) APC, with $K = 12$ (``APC-a1'')~\cite{Tonkoski12}; and, iv) APC with inverter selection (``APC-a2''), where $\lambda_p$ in~\eqref{SGlasso_powers} is chosen such that 7 inverters are controlled. It is interesting to note that the minimum PF constraints were \emph{not} enforced for the RPC strategy; in fact, $(\mathrm{P}4)$ is infeasible for a minimum PF constraint higher than 0.3. Furthermore, the case where $c_\phi = 1$, $b_h = c_\rho = 1$ (with all other coefficients set to $0$), and $\lambda$, $\lambda_p$ chosen such that $7$ inverters are selected is also considered. This represents the situation where
the utility attempts to minimize the overall power losses in the network, namely the sum of the line losses plus the curtailed power. The resultant OID and APC strategies are marked as ``OID-d3'' and ``APC-a3,'' respectively. For a fair comparison with RPC, no minimum PF constraints are enforced in OID-d3. As can be seen from Fig.~\ref{F_voltage}, voltage regulation is effected in all the considered cases, with the flattest voltage profile for APC-a1. However, these strategies yield different active power losses in the network, as well as \emph{overall} active power losses (that is, the power lost in the lines plus the curtailed power). These two quantities are compared in Fig.~\ref{F_power_loss}. When there is no price associated with the active power curtailed, RPC is the one that yields the highest power losses in the network (as hinted in~\cite{Tonkoski11}); however, APC yields the highest power losses on the lines when $b_h = c_\rho = 1$. Interestingly, strategy ``OID-d3''  yields the \emph{lowest} overall active power losses, thus demonstrating the merits of the proposed OID approach. On the other hand, ``OID-d2'' yields a good tradeoff between voltage profile flatness and power loss when $b_h < c_\rho$. Table~\ref{tab:day_results} collects the energy loss in the network (in the column labeled ``Network''), the energy curtailed by the inverters (in the column labeled ``Curtailed'') and the total energy loss (in the column labeled ``Overall'') for the simulated day. The accumulated energy loss in the business-as-usual approach with no ancillary services is also reported for comparison purposes.

\begin{table}[t]
\caption{Energy loss in the network, and curtailed PV energy for the simulated day}
\vspace{-.2cm}
\begin{center}
\begin{tabular}{l||c|c|c}
&  Network & Curtailed & Overall \\
&  [kWh] &  [kWh] & [kWh] \\
\hline \hline
No control & 97.4 & 0 & 97.4 \\
\hline
RPC-r1 & 15.4 & 0 & 15.4 \\
RPC-r2 & 18.0 & 0 & 18.0 \\ 
\hline
($c_\phi = 0$) & & & \\
APC-a1 & 8.1 & 118.3 & 126.4  \\
APC-a2 &  43.4 & 75.4 &  118.8  \\
OID-d1 & 4.1 & 116.6 & 120.7 \\
OID-d2 & 6.5 & 20.0 & 26.5 \\
\hline
($c_\phi = c_\rho) $ & & & \\
APC-a3 & 54.4 & 29.3 &  83.7 \\
OID-d3 & 13.7 & 0.5 & \textbf{14.2} \\
\end{tabular}
\end{center}
\vspace{-.2cm}
\label{tab:day_results}
\end{table}%

\begin{figure}[t]
\centering
  \subfigure[]{\includegraphics[width=0.42\textwidth]{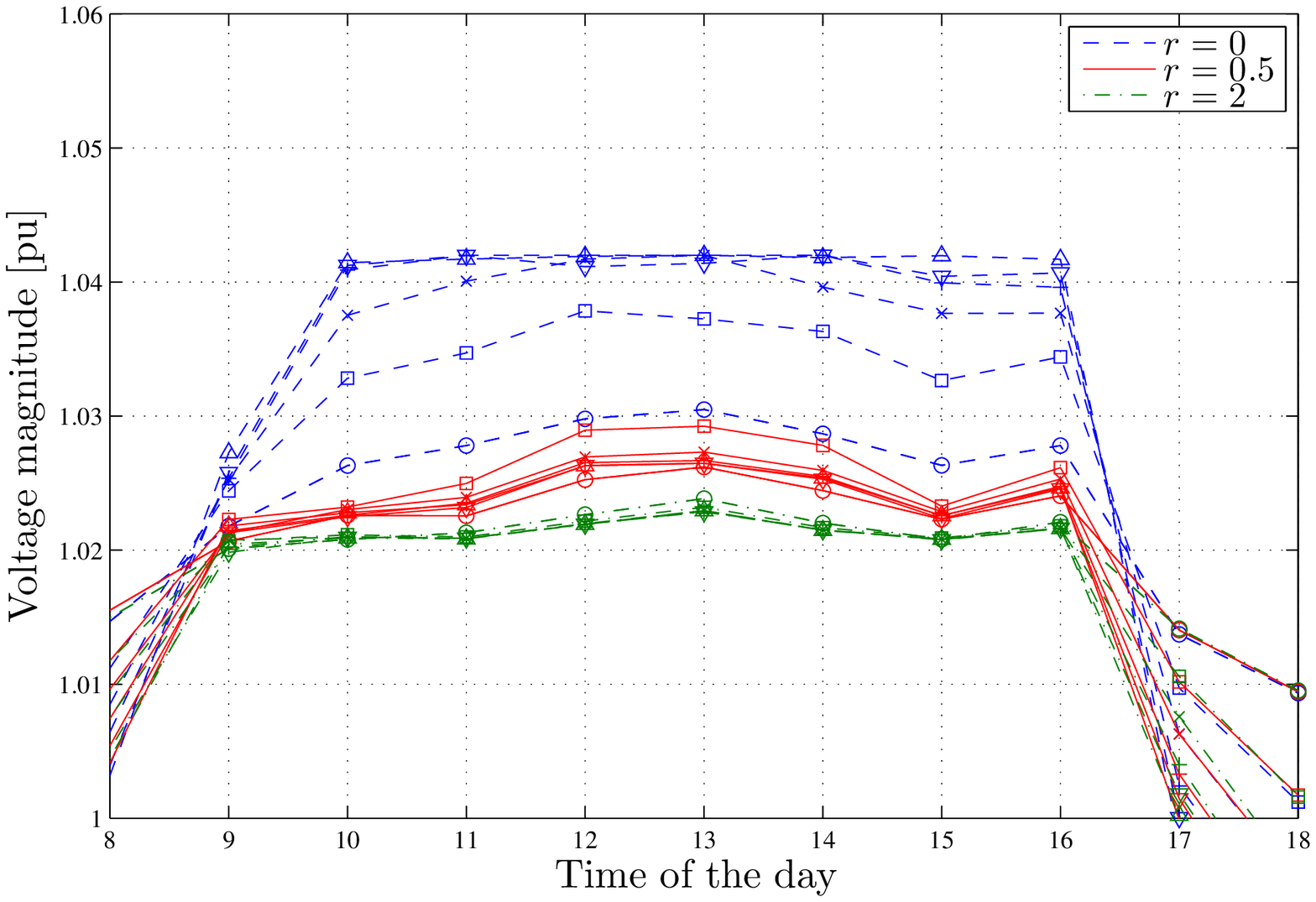}} 
  \subfigure[]{\includegraphics[width=0.405\textwidth]{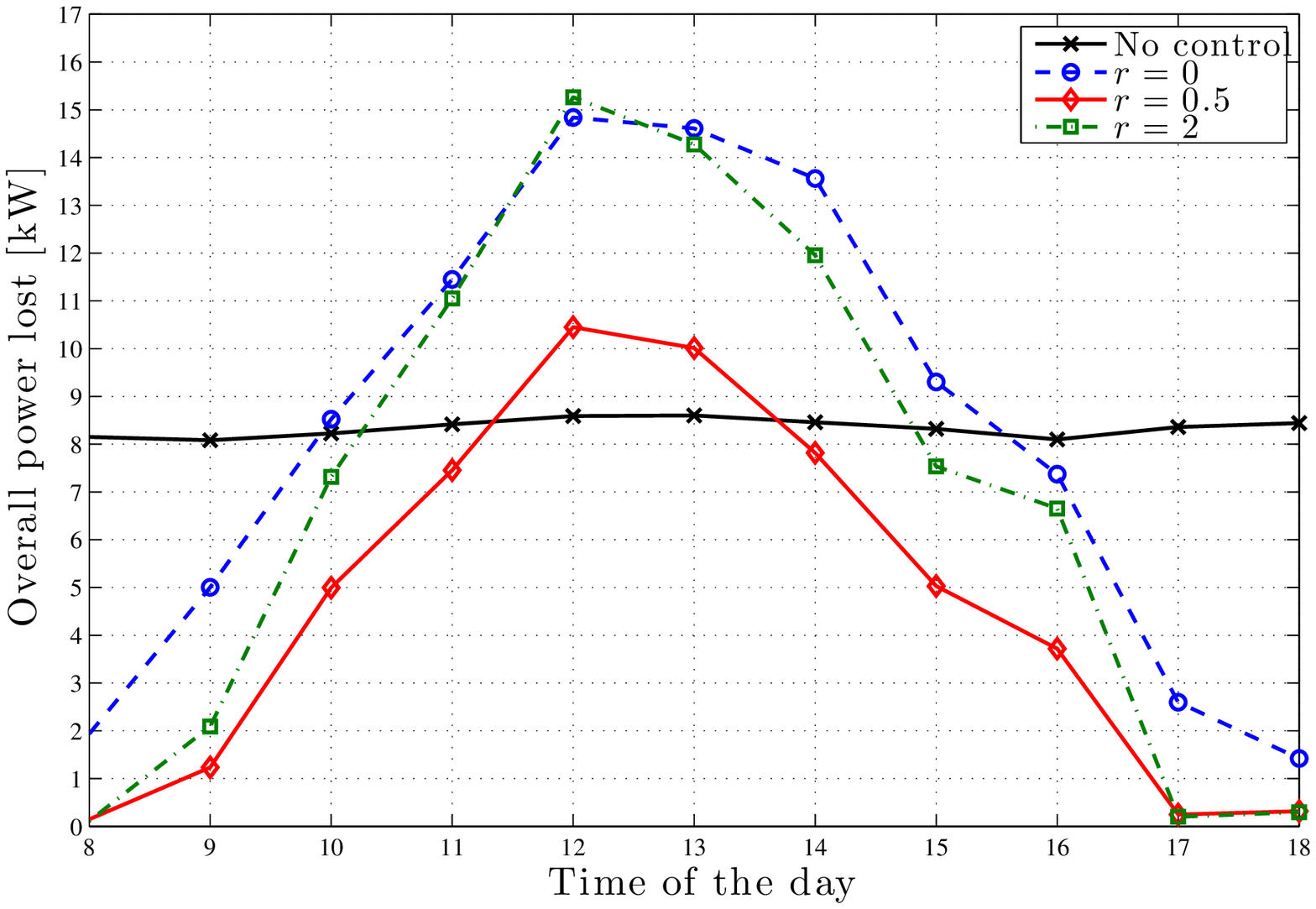}} 
  \caption{Voltage magnitude at $\{\mathrm{H}_i, i = 1,3,5,7,9,11\}$ (a) and overall power lost (b) for different values of $r := c_\nu/c_\rho$. In this setup, $c_\phi = 0$, $\lambda = 0.8$.} 
\label{F_comparison}
\end{figure}

To demonstrate the increased flexibility offered by the function $c_\nu \nu(\bv)$ in~\eqref{eq:CostFunction}, Fig.~\ref{F_comparison} depicts voltage profiles at a few houses during the course of the day for different values of $r := c_\nu/c_\rho$, for $c_\phi = 0, c_\rho = 1$, from which it is evident that deviations from the average can be minimized
by increasing $r$. Fig.~\ref{F_comparison}(b) illustrates that the lowest overall power losses are obtained for $r = 0.5$ (demonstrating that $r$ cannot be increased indiscriminately without considering other optimization objectives).

To better highlight the advantages of the proposed OID method, a long term impact analysis over the course of a year is also performed. 
To this end, the hourly profiles for available solar powers are generated using SAM, based on the TMY data for Minneapolis, MN, and the hourly load profiles available for the whole year in the Open Energy Info database for the Twin Cities, MN.  Table~\ref{tab:annual_results} reports the energy loss in the network and the energy curtailed over the whole year. Two setups are considered: i)  $b_h = c_\phi = 0, c_\rho \neq 0$ (i.e., the active power losses in the network are minimized); and, ii) $b_h  =1, c_\phi = c_\rho$ (i.e., the overall power lost is to be minimized). In the first setup, it can be clearly seen that the OID strategy yields the \emph{lowest} losses in the network. The proposed scheme outperforms the RPC and APC strategies also in the second case, since it yields the \emph{lowest} overall losses. 

To estimate potential economic savings over the course of a year, consider using the ``average retail price for electricity to ultimate customers by end-use sector'' available in the monthly reports of the U.S. Energy Information Administration, for the year 2012, for the state of Minnesota.\footnote{[Online] Available at: \texttt{http://www.eia.gov/totalenergy}} Let $c_{\textrm{retail}}$ denote the average retail price, and consider solving the OID problem over the whole year in the  following two cases: i) $\kappa(\bV, \bp_c, \bq_s) = c_{\textrm{retail}} \trace(\bL \bV)$ 
(that is, only the economic losses in the network are considered); and, ii) $\kappa(\bV, \bp_c, \bq_s) = c_{\textrm{retail}} \trace(\bL \bV) + c_{\textrm{retail}}  \sum_{h \in \cH} P_{c,h}$  (that is, the economic loss emerges from both the losses in the network and the active power curtailed). Table~\ref{tab:annual_results2} summarizes the overall economic losses over the year. It can be clearly seen that, by enabling a joint optimization of the active and reactive power production, OID provides the most economic savings.

\begin{table}[t]
\caption{Energy loss in the network, and curtailed PV energy for the whole year}
\vspace{-.2cm}
\begin{center}
\begin{tabular}{l||c|c|c}
&  Network & Curtailed & Overall \\
& [kWh] & [kWh] & [kWh]  \\
\hline \hline
RPC-r1 & 5040 & 0 & 5040 \\
RPC-r2 & 4710 & 0 & 4710 \\ 
\hline
($b_h = 0 $) & & & \\
APC-a1 & 965 & 34170 & 43820 \\
APC-a2 &  3148 & 23512 &  26660 \\
OID-d1 & 938 & 34700 & 35638 \\
OID-d2 & \textbf{913} & 10862 & 11775 \\
\hline
($c_\phi = c_\rho) $ & & & \\
APC-a3 & 10387 & 3864 &  14251 \\
OID-d3 & 4237 & 143 & \textbf{4380} \\
\end{tabular}
\end{center}
\vspace{-.2cm}
\label{tab:annual_results}
\end{table}%

\begin{table}[t]
\caption{Economic losses over one year [\$]}
\vspace{-.2cm}
\begin{center}
\begin{tabular}{c||r|r|r}
Case & RPC & APC & OID  \\
\hline 
$c_\phi = 0, c_\rho = c_{\textrm{retail}}$ & 540.11 & 355.20 & \textbf{102.94} \\
$b_h = 1, c_\phi = c_\rho = c_{\textrm{retail}}$ & 540.11 &  1520.45 & \textbf{499.46} \\
\end{tabular}
\end{center}
\vspace{-.2cm}
\label{tab:annual_results2}
\end{table}%

\section{Concluding Remarks and Future Directions }
\label{sec:Conclusions}

A framework to facilitate high PV penetration in existing residential
distribution networks is proposed in this paper. Through OID, the inverters that have to be controlled in order to enable voltage regulation are identified, and their optimal active- and reactive-power set points are obtained. The OID problem was cast as an SDP, which is efficiently solvable with a computational time on the order of the sub-second. Overall, the novel OID framework demonstrates that system-level real-time optimization can facilitate the integration of PV systems in existing low-voltage distribution systems, and calls for instituting communication capabilities with residential-scale PV inverters, in order to enable real-time provisioning of ancillary services. Future efforts include investigating the impact of forecast and communication errors on the inverter dispatch task. In addition, the benefits of the proposed OID method can be effectively communicated to practicing engineers by targeting experimental demonstrations in real-world setups.

\appendix

\subsection{Derivation of $(\mathrm{P}3)$}
\label{sec:appA}

An SDP in standard form involves the minimization of a linear function, subject to linear (in)equalities and linear matrix inequalities (LMIs)~\cite{Vandenberghe96}. The SDP reformulation of $(\mathrm{P}2)$ relies on the \emph{Schur complement}. For the symmetric matrix $\bM := \left[\begin{array}{ccc}
\bA & \bB \\
\bB^\mathcal{T} & \bC 
\end{array}
\right]$, where $\bC$ is invertible, the Schur complement of $\bA$ is given by $\bS = \bC-\bB^\mathcal{T}\bA^{-1}\bB$. It follows that $\bM$ is positive definite if and only if $\bS$ and $\bA$ are positive definite (or, if and only if $S\geq 0$ when $A, B, C$ are scalars); this result has been leveraged in the derivation of $\mathrm{(P3)}$ as detailed next.

Aiming for an SDP formulation of $(\mathrm{P}2)$, consider introducing the auxiliary variable $\upsilon \geq 0$, replacing $ \| \bPi \diag(\bV) \|_2 $ with $\upsilon$ in the cost of $(\mathrm{P}3)$, and adding the constraint $ \| \bPi \diag(\bV) \|_2 \leq \upsilon$. Then, upon rewriting this constraint as $\upsilon - \upsilon^{-1} \| \bPi \diag(\bV) \|_2^2 \geq 0$, one can readily obtain~\eqref{SDPvoltage} by using the Schur complement with $A = C = \upsilon$ and $\bB = \bPi \diag(\bV)$. Similarly~\eqref{SDPnorm2} can be obtained by setting $\bA = w_h \bI_{2}$, $\bB = [P_{c,h}, Q_{s,h}]^\cT$, and $C = w_h$; \eqref{SDPConstCostActive} can be obtained by setting $A = -1$, $B = \sqrt{a_h}P_{c,h}$, and $C = b_h P_{c,h} - y_h$; and \eqref{SDPConstReactive} can be obtained by setting $\bA = \bI_{2}$, $\bB = [Q_{s,h}, \bar{P}_{h} - P_{c,h}]^\cT$, and $C = - S_{h}^2$. 

Notice finally that the positive semi-definiteness and rank constraints~\eqref{SDPpositive}--\eqref{SDPRank} jointly ensure that there always exists a vector of voltages $\bv$ such that $\bV = \bv \bv^\cH$ for any feasible $\bV \in \mathbb{C}^{N+1 \times N+1}$ (see e.g.,~\cite{luospmag10}).

\subsection{Soft-thresholding on the inverter set point}
\label{sec:appB}

To rigorously demonstrate the PV-inverter selection capability offered by the proposed relaxed OID problem, results from duality theory~\cite{BertsekasConvexAnOpt} are leveraged next to derive closed-form expressions for the optimal inverter set points. 

Define the real-valued vector $\bx_h := [P_{c,h}, Q_{s,h}]^\cT$, and consider rewriting~\eqref{eq:phiP} as $c_\phi \phi(\bp_c, \bq_s) = \sum_{h \in \cH} \bx_h^\cT \bA_h \bx_h + \bb_h^\cT \bx_h$, with $\bA_h := \diag([c_\phi a_h, 0])$ and $\bb_h := [c_\phi b_h,0]^\cT$. Notice further that constraint $Q_{s,h}^2 + (\bar{P}_{h} - P_{c,h})^2 - S_{h}^2 \leq 0$ [cf.~\eqref{mg-PVq2}] can be re-expressed in quadratic form as  
$\bx_h^\cT \bx_h + \bd_h^\cT \bx_h + (\bar{P}_{h}^2 - S_h^2) \leq 0$, with $\bd_h := [-2 \bar{P}_h , 0]^\cT$.  Suppose for simplicity that $c_\nu = 0$ and no PF constraints are imposed, and consider the following relaxed OID problem:
\begin{subequations}
\begin{align} 
& \hspace{-.3cm} \mathrm{(P5)}  \min_{\{\bx_h\} , \bV} c_\rho \trace(\bL \bV) + \sum_{h \in \cH} \left( \bx_h^\cT \bA_h \bx_h + \bb_h^\cT \bx_h  + \lambda \|\bx_2\|_2 \right)  \nonumber \\
& \mathrm{subject\,to}~\eqref{SDPconst}, \bV \succeq \mathbf{0},~\mathrm{and} \nonumber \\
& \hspace{1.2cm} 0 \leq P_{c,h}  \leq  \bar{P}_{h} \label{mg4-PVp2} \\
& \hspace{1.2cm} \bx_h^\cT \bx_h  + \bd_h^\cH \bx_h + (\bar{P}_{h}^2 - S_h^2) \leq 0 .  \label{mg4-PVq} 
\end{align}
\end{subequations}
Problems $(\mathrm{P}4)$ and $(\mathrm{P}5)$ are equivalent, and their globally optimal solution coincide~\cite{BertsekasConvexAnOpt}. Further, it can be shown that Slater's condition holds~\cite{LavaeiLow}, and thus $(\mathrm{P}5)$ has zero duality gap~\cite[Ch.~6]{BertsekasConvexAnOpt}. Let $\{\mu_h\}, \{\bar{\mu}_h\}$ denote the multipliers associated with~\eqref{SDPconstP} and~\eqref{SDPconstQ}, respectively; $\{\varphi_h, \bar{\varphi}_h\}$ the ones with~\eqref{SDPconstV}; and, $\{\gamma_h\}$, $\{\nu_h, \bar{\nu}_h\}$ the ones with the inverter-related constraints~\eqref{mg4-PVp2} and~\eqref{mg4-PVq}, respectively. Further, let $\cL(\bV, \{\bx_h\}, \{\mu_h, \bar{\mu}_h, \varphi_h, \bar{\varphi}_h, \gamma_h, \nu_h, \bar{\nu}_h\})$ denote the Lagrangian of $(\mathrm{P}5)$.         

With the optimal dual variables $\{\mu_h^\textrm{opt}, \bar{\mu}_h^\textrm{opt}, \varphi_h^\textrm{opt}, \bar{\varphi}_h^\textrm{opt},$ $\gamma_h^\textrm{opt}, \nu_h^\textrm{opt}, \bar{\nu}_h^\textrm{opt}\}$ known, the Lagrangian optimality condition~\cite[Prop.~6.2.5]{BertsekasConvexAnOpt} asserts that $\bV^{\textrm{opt}}, \{\bx_h^{\textrm{opt}}\}$ can be found as $\min_{\{\bx_h\} , \bV \succeq \mathbf{0}} \cL(\bV, \{\bx_h\}, \{\mu_h^\textrm{opt}, \bar{\mu}_h^\textrm{opt}, \varphi_h^\textrm{opt}, \bar{\varphi}_h^\textrm{opt}, \gamma_h^\textrm{opt}, \nu_h^\textrm{opt}, \bar{\nu}_h^\textrm{opt}\})$. Thus, exploiting the decomposability of the Lagrangian, it turns out that the optimal setpoint $\bx_h$ for inverter $h$ is given as the solution of the sub-problem:
\begin{align}
\label{eq:GLasso}
\hspace{-.2cm} \bx_{h}^{opt} = \arg \min_{\bx} \frac{1}{2} \bx^\cT \bQ_h \bx + \bc_h^\cT \bx + \lambda \|\bx\|_2
\end{align}
where $\bQ_h := 2 \bA_h + 2 \gamma_h^{\textrm{opt}} \bI_{2}$ and $\bc_h := \bb_h + \bmu_h + \bnu_h + \gamma_h^{\textrm{opt}} \bd_h$, with $\bmu_h := [\mu_h^{\textrm{opt}}, - \bar{\mu}_h^{\textrm{opt}}]^\cT$ and $\bnu_h := [\bar{\nu}_h^{\textrm{opt}} - \nu_h^{\textrm{opt}}, 0]^\cT$. Notice that $\bQ_h \succeq \mathbf{0}$. Then, from~\cite[Thm.~1]{Wiesel11}, it follows that the optimal set points $\bx_h$ are given by 
 the following shrinkage and thresholding vector operation   
\begin{align}
 \bx_{h}^{\textrm{opt}}  & =  \mathbb{I}_{\left\{\|\bc_h\|_2 > \lambda \right\}}
\left[ -\eta \left(\eta \bQ_h  + \frac{\lambda^2}{2} \bI_{2} \right)^{-1} \bc_h \right] \label{Optimal_pq}   
\end{align}
with $\mathbb{I}_{\left\{ A \right\}} = 1$ if event $A$ is true and zero otherwise, and $\eta \in \mathbb{R}^+$  the solution of the scalar optimization problem:
\begin{align}
 \min_{\eta \geq 0} \, \eta - \frac{\eta}{2} \bc_h^\cT \left(\eta \bQ_h  + \frac{\lambda^2}{2} \bI_{2} \right)^{-1} \bc_h \, .
 \label{Optimal_eta}  
 \end{align}
It can be clearly deduced that $\bx_{h}^{\textrm{opt}} = \mathbf{0}$ when $\lambda \geq \|\bc_h\|_2$; that is, inverter $h$ operates at the unitary-PF set point $(\bar{P}_h, 0)$.  On the other hand, when  $[P_{c,h}, Q_{s,h}]^\mathcal{T} \neq \mathbf{0}$, the operating point of inverter $h$ is given by $(\bar{P}_h - P_{c,h}, Q_{s,h})$. Equation~\eqref{Optimal_pq} also explains why, with the decreasing of $\lambda$, the number of controlled inverters increases.

\bibliographystyle{IEEEtran}
\bibliography{biblio}

\begin{thebibliography}{10}
\providecommand{\url}[1]{#1}
\csname url@samestyle\endcsname
\providecommand{\newblock}{\relax}
\providecommand{\bibinfo}[2]{#2}
\providecommand{\BIBentrySTDinterwordspacing}{\spaceskip=0pt\relax}
\providecommand{\BIBentryALTinterwordstretchfactor}{4}
\providecommand{\BIBentryALTinterwordspacing}{\spaceskip=\fontdimen2\font plus
\BIBentryALTinterwordstretchfactor\fontdimen3\font minus
  \fontdimen4\font\relax}
\providecommand{\BIBforeignlanguage}[2]{{%
\expandafter\ifx\csname l@#1\endcsname\relax
\typeout{** WARNING: IEEEtran.bst: No hyphenation pattern has been}%
\typeout{** loaded for the language `#1'. Using the pattern for}%
\typeout{** the default language instead.}%
\else
\language=\csname l@#1\endcsname
\fi
#2}}
\providecommand{\BIBdecl}{\relax}
\BIBdecl

\bibitem{Sherwood12}
L.~Sherwood, ``{U.S.} solar market trends 2012,'' Jul. 2013, [Online] Available
  at \texttt{http://www.irecusa.org}.

\bibitem{Liu08}
E.~Liu and J.~Bebic, ``Distribution system voltage performance analysis for
  high-penetration photovoltaics,'' Feb. 2008, {NREL} Technical Monitor: B.
  Kroposki. Subcontract Report {NREL/SR}-581-42298.

\bibitem{Tonkoski12}
R.~Tonkoski, R.~Turcotte, and T.~H.~M. El-Fouly, ``Impact of high {PV}
  penetration on voltage profiles in residential neighborhoods,'' \emph{IEEE
  Trans. on Sust. Energy}, vol.~3, no.~3, pp. 518--527, Jul. 2012.

\bibitem{Turitsyn11}
K.~Turitsyn, P.~Sulc, S.~Backhaus, and M.~Chertkov, ``Options for control of
  reactive power by distributed photovoltaic generators,'' \emph{Proc. of the
  IEEE}, vol.~99, no.~6, pp. 1063--1073, 2011.

\bibitem{Standard1547}
``{IEEE} 1547 standard for interconnecting distributed resources with electric
  power systems,'' [Online]. Available:
  \url{http://grouper.ieee.org/groups/scc21/dr_shared/}.

\bibitem{Cpuc12}
{California Public Utilities Commission}, ``Advanced inverter technologies
  report,'' Jan. 2013, [Online] Available at \texttt{http://www.cpuc.ca.gov}.

\bibitem{Nerc12}
{North American Electric Reliability Corporation}, ``Special reliability
  assessment: Interconnection requirements for variable generation,'' Mar.
  2012, [Online] Available at \texttt{http://www.nerc.com}.

\bibitem{Carvalho08}
P.~Carvalho, P.~Correia, and L.~Ferreira, ``Distributed reactive power
  generation control for voltage rise mitigation in distribution networks,''
  \emph{IEEE Trans. Power Syst.}, vol.~23, no.~2, pp. 766--772, May 2008.

\bibitem{DeBrabandere04}
M.~A. Mahmud, M.~J. Hossain, H.~R. Pota, and A.~B.~M. Nasiruzzaman, ``Voltage
  control of distribution networks with distributed generation using reactive
  power compensation,'' in \emph{Conf. on IEEE Ind. Electr. Soc.}, 2011, pp.
  985--990.

\bibitem{Cagnano11}
A.~Cagnano, E.~D. Tuglie, M.~Liserre, and R.~A. Mastromauro, ``Online optimal
  reactive power control strategy of {PV} inverters,'' \emph{IEEE Trans. on
  Ind. Electron.}, vol.~58, no.~10, pp. 4549--4558, Oct. 2011.

\bibitem{Farivar12}
M.~Farivar, R.~Neal, C.~Clarke, and S.~Low, ``Optimal inverter {VAR} control in
  distribution systems with high {PV} penetration,'' in \emph{IEEE PES General
  Meeting}, San Diego, CA, Jul. 2012.

\bibitem{Aliprantis13}
P.~Jahangiri and D.~C. Aliprantis, ``Distributed {Volt/VAr} control by {PV}
  inverters,'' \emph{IEEE Trans. Power Syst.}, vol.~28, no.~3, pp. 3429--3439,
  Aug. 2013.

\bibitem{Tonkoski11}
R.~Tonkoski, L.~A.~C. Lopes, and T.~H.~M. El-Fouly, ``Coordinated active power
  curtailment of grid connected {PV} inverters for overvoltage prevention,''
  \emph{IEEE Trans. on Sust. Energy}, vol.~2, no.~2, pp. 139--147, Apr. 2011.

\bibitem{Tonkoski11Renewable}
R.~Tonkoski and L.~A.~C. Lopes, ``Impact of active power curtailment on
  overvoltage prevention and energy production of {PV} inverters connected to
  low voltage residential feeders,'' \emph{Renewable Energy}, vol.~36, no.~12,
  pp. 3566--3574, Dec. 2011.

\bibitem{ADG-TPWRS12}
B.~A. Robbins, C.~N. Hadjicostis, and A.~D. Dom\'{i}nguez-Garc\'{i}a, ``A
  two-stage distributed architecture for voltage control in power distribution
  systems,'' \emph{IEEE Trans. Power Syst.}, vol.~28, no.~2, pp. 1470--1482,
  2012.

\bibitem{YuLi06}
M.~Yuan and Y.~Lin, ``Model selection and estimation in regression with grouped
  variables,'' \emph{J. of the Royal Stat. Soc.}, vol.~68, pp. 49--67, 2006.

\bibitem{Wiesel11}
A.~T. Puig, A.~Wiesel, G.~Fleury, and A.~O. Hero, ``Multidimensional
  shrinkage-thresholding operator and group {LASSO} penalties,'' \emph{IEEE
  Sig. Proc. Letters}, vol.~18, no.~6, pp. 363--366, Jun. 2011.

\bibitem{Bai08}
X.~Bai, H.~Wei, K.~Fujisawa, and Y.~Wang, ``Semidefinite programming for
  optimal power flow problems,'' \emph{Int. J. of Electrical Power \& Energy
  Systems}, vol.~30, no. 6--7, pp. 383--392, Jul.-Sep. 2008.

\bibitem{LavaeiLow}
J.~Lavaei and S.~H. Low, ``Zero duality gap in optimal power flow problem,''
  \emph{IEEE Trans. Power Syst.}, vol.~1, no.~1, pp. 92--107, Feb. 2012.

\bibitem{Tse12}
A.~Y. Lam, B.~Zhang, A.~Dom\'{i}nguez-Garc\'{i}a, and D.~Tse, ``Optimal
  distributed voltage regulation in power distribution networks,'' 2012,
  [Online] Available at \texttt{http://arxiv.org/abs/1204.5226v1}.

\bibitem{Dallanese-TSG13}
E.~Dall'Anese, H.~Zhu, and G.~B. Giannakis, ``Distributed optimal power flow
  for smart microgrids,'' \emph{IEEE Trans. Smart Grid}, vol.~4, no.~3, pp.
  1464--1475, Sep. 2013.

\bibitem{luospmag10}
Z.-Q. Luo, W.-K. Ma, A.~M.-C. So, Y.~Ye, and S.~Zhang, ``Semidefinite
  relaxation of quadratic optimization problems,'' \emph{IEEE Sig. Proc. Mag.},
  vol.~27, no.~3, pp. 20--34, May 2010.

\bibitem{Lavaei_tree}
J.~Lavaei, D.~Tse, and B.~Zhang, ``Geometry of power flows and optimization in
  distribution networks,'' in \emph{IEEE PES General Meeting.}, San Diego, CA,
  2012.

\bibitem{Vandenberghe96}
L.~Vandenberghe and S.~Boyd, ``Semidefinite programming,'' \emph{{SIAM}
  Review}, vol.~38, no.~1, pp. 49--95, Mar. 1996.

\bibitem{PaudyalyISGT}
S.~Paudyaly, C.~A. Canizares, and K.~Bhattacharya, ``Three-phase distribution
  {OPF} in smart grids: Optimality versus computational burden,'' in \emph{2nd
  IEEE PES Intl. Conf. and Exhibition on Innovative Smart Grid Technologies},
  Manchester, UK, Dec. 2011.

\bibitem{Key12}
T.~Key and B.~Seal, ``{Inverters to Provide Grid Support (DG and Storage)},''
  in \emph{Proc. of of 5th International Conference on Integration of Renewable
  and Distributed Energy Resources}, Dec. 2012.

\bibitem{Bolognani13}
S.~Bolognani and S.~Zampieri, ``A distributed control strategy for reactive
  power compensation in smart microgrids,'' \emph{IEEE Trans. on Autom.
  Control}, vol.~58, no.~11, pp. 2818--2833, Nov. 2013.

\bibitem{GSO08}
J.~D. Glover, M.~S. Sarma, and T.~J. Overbye, \emph{Power System Analysis and
  Design}.\hskip 1em plus 0.5em minus 0.4em\relax Thomson Learning, 2008.

\bibitem{Braun10}
M.~Braun, J.~K\"unschner, T.~Stetz, and B.~Engel, ``Cost optimal sizing of
  photovoltaic inverters -- influence of of new grid codes and cost
  reductions,'' in \emph{Proc. of 25th Europ. {PV} Solar Energy Conf. and
  Exhib.}, Valencia, Spain, Sep. 2010.

\bibitem{SolarBridge}
\BIBentryALTinterwordspacing
{SolarBridge Technologies}. (Retrieved, October 2013) Solarbridge power
  manager. [Online]. Available:
  \url{http://solarbridgetech.com/products/our-solution/solarbridge-power-mana%
ger/}
\BIBentrySTDinterwordspacing

\bibitem{Enphase}
\BIBentryALTinterwordspacing
{Enphase Energy}. (Retrieved, October 2013) Envoy communications gateway.
  [Online]. Available: \url{http://enphase.com/products/envoy/}
\BIBentrySTDinterwordspacing

\bibitem{Farivar13}
M.~Farivar and S.~H. Low, ``Branch flow model: Relaxations and convexification
  {(Part I)},'' \emph{IEEE Trans. Power Syst.}, vol.~28, no.~3, pp. 2554--2564,
  2013.

\bibitem{Tib96}
R.~Tibshirani, ``Regression shrinkage and selection via the {L}asso,'' \emph{J.
  Royal Stat. Soc.}, vol.~58, no.~1, pp. 267--288, 1996.

\bibitem{King(SANDIA):2007}
D.~L. King, S.~Gonzalez, G.~M. Galbraith, and W.~E. Boyson, ``Performance model
  for grid-connected photovoltaic inverters,'' Sandia National Laboratories,
  Tech. Rep., September 2007.

\bibitem{Cady2012}
S.~T. Cady, D.~Mestas, and C.~Cirone, ``Engineering systems in the rehome: A
  net-zero, solar-powered house for the {U.S.} {D}epartment of {E}nergy's 2011
  {S}olar {D}ecathlon,'' in \emph{{IEEE} Power and Energy Conf. at Illinois},
  Univ. Illinois at Urbana-Champaign, 2012.

\bibitem{Dhople2009}
S.~V. Dhople, J.~L. Ehlmann, C.~J. Murray, S.~T. Cady, and P.~L. Chapman,
  ``Engineering systems in the gable home: A passive, net-zero, solar-powered
  house for the {U.S.} {D}epartment of {E}nergy's 2009 {S}olar {D}ecathlon,''
  in \emph{{IEEE} Power and Energy Conf. at Illinois}, Univ. Illinois at
  Urbana-Champaign, 2010.

\bibitem{BertsekasConvexAnOpt}
D.~P. Bertsekas, A.~Nedi\'{c}, and A.~Ozdaglar, \emph{Convex Analysis and
  Optimization}.\hskip 1em plus 0.5em minus 0.4em\relax Athena Scientific,
  2003.

\end{thebibliography}

\end{document}